\documentclass{svmult}
\usepackage{epic}
\usepackage{eepic}
\usepackage{epsf} 
\usepackage{amscd}
\usepackage{amsmath}
\usepackage{epsfig}
\usepackage{graphicx}
\usepackage{amsbsy}
\usepackage{amssymb}
\usepackage{hyperref}
\usepackage{mathrsfs}
\usepackage{mathtools}
\usepackage{euscript}
\usepackage{eufrak}
\usepackage{upgreek}

\usepackage{mathptmx}
\usepackage{helvet}
\usepackage{courier}
\usepackage{makeidx}
\usepackage{multicol}
\usepackage{footmisc}

\makeatother

\usepackage[english]{babel}
\newcommand{\uhp}{\mathbb{H}}
\newcommand{\SSL}{\boldsymbol{\Gamma}}

\newcommand{\mr}[4]{\rho\!\left(\begin{array}{cc}
 #1  &  #2\\
#3  &  #4\end{array} \right)}

\newcommand{\FA}{\mathbb{X}\!}
\newcommand{\Fa}{\mathbb{X}}
\newcommand{\FB}{{\Lambda}}
\newcommand{\FC}{\mathrm{z}}
\newcommand{\FE}{\mathcal{P}}
\newcommand{\DO}{\nabla}
\newcommand{\Der}{\mathrm{D}}
\newcommand{\varmu}{\upsilon}

\newcommand{\bbH}{\mathbb{H}}
\newcommand{\bbQ}{\mathbb{Q}}
\newcommand{\bbC}{\mathbb{C}}
\newcommand{\bbK}{\mathbb{K}}
\newcommand{\bbP}{\mathbb{P}}

\newcommand{\bbY}{\mathbb{Y}}
\newcommand{\bbX}{\mathbb{X}}
\newcommand{\bbZ}{\mathbb{Z}}

\newcommand{\Ii}{\mathrm{i}}
\newcommand{\cA}{\mathcal{A}}

\newcommand{\cK}{\mathcal{K}}

\newcommand{\cF}{\mathcal{F}}

\newcommand{\cM}{\mathcal{M}}
\newcommand{\cP}{\mathcal{P}}

\newcommand{\cT}{\mathcal{T}}

\newcommand{\bG}{{{\Gamma}}}

\newcommand{\bF}{{\mathrm{F}}}
\newcommand{\sdsum}{\ni\!\!\!\!\!\!\!\vrule height5pt depth0pt width0.4pt\,\,\,\,\,\,}

\newcommand{\map}[3]{#1\!:#2\!\rightarrow\!#3}
\usepackage[utopia]{mathdesign}
\DeclareMathSymbol{\coprodd}{\mathop}{largesymbols}{"60}
\def\frakm{\mathrm{m}}

\begin{document}
\title*{The theory of vector-valued modular forms for the modular group}
\author{Terry  Gannon }
\institute{Terry Gannon \at Department of Mathematics, University of Alberta,
 Edmonton, Alberta, Canada T6G 2G1, \email{tgannon@math.ualberta.ca}}
\titlerunning{Vector-valued modular forms}
\maketitle

\abstract{We explain the basic ideas, describe with proofs the
main results, and demonstrate the effectiveness,
of an evolving theory of vector-valued modular forms (vvmf).  To keep the exposition
concrete, we restrict here to the special case of the modular group. Among other things,
we construct vvmf for arbitrary multipliers, solve the Mittag-Leffler problem here, establish
Serre duality and find a dimension formula for holomorphic vvmf, all in far greater generality than
has been done elsewhere. More important, the new ideas involved are sufficiently simple
and robust that this entire theory extends directly to any genus-0 Fuchsian group.} 
%\end{abstract}

\section{Introduction}

Even the most classical modular forms (e.g. the Dedekind eta $\eta(\tau)$)
need a multiplier, but this multiplier is typically a number (i.e. a 
1-dimensional projective representation of some discrete group like
$\Gamma={\rm SL}_2(\bbZ)$).  Simple examples of \textit{vector-valued}
modular forms (vvmf) for SL$_2(\bbZ)$ are the weight-$\frac{1}{2}$ Jacobi theta functions
$\Theta(\tau)=(\theta_2(\tau),\theta_3(\tau),\theta_4(\tau))^t$, which obey for instance
\begin{equation}\label{ex1}
\Theta(-1/\tau)
=\sqrt{\frac{\tau}{\Ii}}\left(\begin{matrix}0&0&1\\ 0&1&0\\ 1&0&0\end{matrix}\right)\Theta(\tau)\ ,\end{equation}
 and $P(\tau)=(\tau,1)^t$, which has weight $-1$ and obeys for instance
 \begin{equation}P(-1/\tau) =\tau^{-1}\left(\begin{matrix}0& -1\\ 1&0\end{matrix}
\right)P(\tau)\ .\label{ex2}\end{equation}

Back in the 1960s Selberg \cite{bib:Se} called for the development of the theory of vvmf, as a way 
to study growth of coefficients of (scalar) modular forms for noncongruence
groups. Since then, the
relevance of vvmf has grown significantly, thanks largely
to the work of Borcherds (see e.g. \cite{bib:Borlift}) and, in physics, the rise of  rational conformal field theory
(RCFT).

In particular, the characters of rational and logarithmic conformal
field theories,  or $C_2$-cofinite vertex operator algebras,  form a weight-0 vvmf for $\bG$
\cite{bib:Zh,bib:Miy}. Less appreciated is that the 4-point functions (conformal blocks) on the 
sphere in RCFT
can naturally be interpreted as vvmf for $\Gamma(2)$,
through the identification of the moduli space of 4-punctured spheres
with $\Gamma(2)\backslash\uhp$. Moreover,
the 1- and 2-point functions on a torus are naturally identified
with vector-valued modular resp. Jacobi forms for $\Gamma$. 
Also, if the RCFT has additional structure, e.g.
$N=2$ supersymmetry or a Lie algebra  symmetry, then the 1-point functions
on the torus can be augmented, becoming vector-valued Jacobi (hence
matrix-valued modular) forms for $\Gamma_0(2)$ or $\Gamma$ \cite{bib:Miy1,bib:KrMa}.

The impact in recent years of RCFT on mathematics makes it difficult to dismiss
 these as esoteric exotica. For example, 1-point torus functions for the orbifold of 
the Moonshine
module $V^\natural$ by subgroups of the Monster contain as very special cases the
Norton series of generalized Moonshine, so a study of them could lead to 
extensions of the Monstrous Moonshine conjectures. For these less
well-known applications to
RCFT, the multiplier $\rho$ will typically have infinite image, and
the weights can be arbitrary rational numbers. Thus the typical classical
assumptions that the weight be half-integral, and the modular form be fixed
by some finite-index subgroup of $\Gamma$, is violated by a plethora
of potentially interesting examples. Hence in the following we do not make those classical
assumptions (nor are they needed). In fact, there should be similar
applications to sufficiently nice non-rational CFT, such as Liouville
theory, where the weights $w$ can be irrational.

In spite of its relevance, the general theory of vvmf has been slow in coming.
Some effort has been made (c.f. \cite{bib:Sk,bib:KM}) to 
 lift to vvmf, classical results like dimension formulas and the `elementary' 
 growth estimates of  Fourier coefficients. Moreover,
differential equations have been recognised as valuable tools for studying vvmf, for
many years (c.f. \cite{bib:AM,bib:MMS,bib:Mil,bib:Mas} to name a few). 

Now, an elementary observation is that  a vvmf  $\bbX(\tau)$ for 
a finite-index subroup $\mathrm{G}$ of $\bG$ can be lifted to one of $\bG$, by inducing the
multiplier. This increases the rank of the vvmf by a factor equal to the index.
This isomorphism tells us that developing a theory of vvmf for $\bG$ gives for free that
of any finite-index subgroup. But more important perhaps, it also shows that the theory of vvmf for $\bG$ contains as a small subclass
the scalar modular forms for noncongruence subgroups. This means that one can only be so successful in lifting results from
the classical ($=$scalar) theory to vvmf. We should be looking for new ideas!

Our approach is somewhat different, and starts from the
heuristic that a vvmf for a Fuchsian group $\bG$ is a lift to $\uhp$ 
of a meromorphic 
section of a flat holomorphic vector
bundle over the (singular) curve $\bG\backslash \uhp$ compactified if necessary by adjoining 
the cusps.  The cusps mean we are not in the
world of algebraic stacks. In place of the order-$d$ ODE on $\bbH$ studied by other authors,
we consider a first order Fuchsian DE on the sphere. \textit{Fuchsian differential equations on
compact curves, and vvmf for Fuchsian groups, are two sides of the same coin.} Another crucial 
ingredient of our theory is the behaviour at the elliptic fixed-points. This has been largely ignored 
in the literature.  For simplicity, this paper restricts to the most familiar (and  important)  case:
$\bG=\mathrm{SL}_2(\bbZ)$, where $\bG\backslash \uhp^*$ is a sphere with 3 conical
singularities (2 at elliptic points and 1 at the cusp). The theory for other Fuchsian  groups is developed elsewhere \cite{bib:BG3,bib:Ga}. 

There are two aspects to the theory: \textit{holomorphic} (no poles anywhere) and \textit{weakly holomorphic} (poles are allowed at the cusps, but only there). 
We address both. We start with weakly holomorphic not because it is more interesting,
but because it is easier, and this makes it more fundamental. There is nothing particularly
special about the cusps from this perspective --- the poles could be allowed at any finitely
many $\Gamma$-orbits, and the theory would be the same. 

Section 3 is the heart of this paper. There we establish existence, using R\"ohrl's solution to
the Riemann--Hilbert problem. We obtain analogues of 
the Birkhoff--Grothendieck  and Riemann--Roch Theorems. We find a dimension formula for
holomorphic vvmf, and are able to quantify the failure of exactness of the functor assigning
to multipliers $\rho$, spaces of holomorphic vvmf. Our arguments are simpler and much more
general than others in the literature. In Section 4 we give several illustrations
of the effectiveness of the theory.

\section{Elementary Remarks}

\subsection{The geometry of the modular group}

Fix $\xi_n=\exp(2\pi\Ii/n)$. Complex powers $z^{w}$ throughout the
paper are defined by $z^{w}=|z|^{w}e^{w\Ii\mathrm{Arg}(z)}$ for $-\pi\le
\mathrm{Arg}(z)<\pi$. We write $\bbC[x]$ for the space of polynomials, $\bbC[[x]]$ for
power series $\sum_{n=0}^\infty a_nx^n$, and $\bbC[x^{-1},x]]$ for Laurent expansions
$\sum_{n=-N}^\infty a_nx^n$ for any $N$.

This paper restricts to the modular group  $\Gamma:=\mathrm{SL}_2(\bbZ)$. Write $\overline{\bG}=\mathrm{PSL}_2(\bbZ)$.
Throughout this paper we use
\begin{equation}S=\mr{0}{-1}{1}{-1}\,,\ \ T=\mr{1}{1}{0}{1}\,,\ \ U=\mr{0}{-1}{1}{-1}=ST^{-1}\,.\end{equation}
Write $\bbH^*$ for the \textit{extended half-plane} $\bbH\cup\bbQ\cup\{\infty\}$. Then
$\Gamma\backslash\bbH^*$ is topologically a sphere. As it is genus 0, it is uniformised by
a Hauptmodul, which can be chosen to be
\begin{equation}\label{Jfunc}J(\tau)=q^{-1}+744+196\,884\,q+21\,493\,760\,q^2+864\,299\,970\,q^3+\cdots\,,\end{equation}
where as always $q=e^{2\pi\Ii\tau}$. For $k\ge 2$ write $E_k(\tau)$ for the Eisenstein series,
normalised so that $E_k(\tau)=1+\ldots$. For $\bG$, the ring of 
{\it weakly holomorphic} (i.e. poles allowed only
at the cusp) modular functions is 
$\bbC[J]$ and the ring of \textit{holomorphic} (i.e. holomorphic everywhere including
the cusp) modular forms is 
$\frakm=\bbC[E_4,E_6]$.

In the differential structure induced by that of $\uhp$,  $\bG\backslash\uhp^*$ will have
a singularity for every orbit of $\overline{\Gamma}$-fixed-points. 
$J(\tau)$ smooths out these 3 singularities. The important $q$-expansion is the local expansion at
one of those singulatities, but it is a mistake to completely ignore the other two, at $\tau=\Ii$ and
$\tau=\xi_6$. 
These elliptic point expansions have been used for at least a century, even if they are largely ignored today.
They play a crucial role in our analysis.

In particular, define $\tau_2=\epsilon_2\,(\tau-\Ii)/(\tau+\Ii)$ and
$\mathrm{j}_2(w;\tau)=E_4(\Ii)^{-w/4}(1-\tau_2/\epsilon_2)^{-w}$, where  $\epsilon_2=\pi \sqrt{E_4(\Ii)}$ and $E_4(\Ii)=3\Gamma(1/4)^8/(2\pi)^6$, and 
define $\tau_3=\epsilon_3\,(\tau-\xi_6)/(\tau-\xi_6^5)$ and
$\mathrm{j}_3(w;\tau)=E_6(\xi_6)^{-w/6}(1-\tau_3/\epsilon_3)^{-w}$, where  $\epsilon_3=\pi\sqrt{3}{E_6(\xi_6)}^{1/3} $ and $E_6(\xi_6)=27\Gamma(1/3)^{18}/(2^9\pi^{12})$. The transformations $\tau\mapsto \tau_i$ for $i=2,3$
map $\uhp$ onto the discs $|\tau_i|<\epsilon_i$, and send $\Ii$ and $\xi_6$ respectively to 0; $\mathrm{j}_2$ and
$\mathrm{j}_3$ are proportional to the corresponding multipliers for a weight-$w$ modular form. 
The elliptic fixed-point $\tau=\Ii$ is fixed by $\left({0\atop 1}{-1\atop 0}\right)$, which sends
$\tau_2\mapsto -\tau_2$. The elliptic fixed-point $\tau=\xi_6$ is
fixed by $\left({0\atop 1}{-1\atop-1}\right)$; it sends $\tau_3\mapsto
\xi_3\tau_3$. If $f$ is a scalar modular form of weight $k$, and we write $\widetilde{f}_i(\tau)
=\mathrm{j}_i(k;\tau)f(\tau_i)$ for $i=2,3$, we get $\widetilde{f}_2(-\tau_2)
=\Ii^{k}\widetilde{f}(\tau_2)$ and $\widetilde{f}_3(\xi_3\tau_3)=\xi_6^k\widetilde{f}_3(\tau_3)$.
We generalise this in Lemma 2.1 below.

We have 
rescaled $\tau_2,\tau_3$ by $\epsilon_2,\epsilon_3$ to clean up the expansions, but this
isn't important in the following. For instance, we have the rational expansions
\begin{align}  
\mathrm{j}_2(4;\tau){E_4}(\tau)&=1+10\tau_2^2/9+{5}\tau_2^4/27+4\tau_2^6/81+19\tau_2^8/5103\cdots\,,\\
\mathrm{j}_2(6;\tau){E_6}(\tau)&=2\tau_2+28\tau_2^3/27+56\tau_2^5/135+28\tau_2^7/405+\cdots\,,\\
J(\tau)&=1728+6912\tau^2_2+11776\tau_2^4+1594112\tau_2^6/135+\cdots
\,.\end{align}
Curiously, the expansion coefficients for $J(\tau)$ at $\Ii$ are
 all \textit{positive} rationals, but infinitely many distinct primes divide 
the denominators. However, 
these denominators arise because of an $n!$ that appears implicitly in these coefficients (see Proposition 17 in \cite{bib:Za}).
Factoring that off, the sequence becomes: 
\begin{align}\nonumber 1728, 10368, &158976, 3586752, 107057664, 4097780928, 193171879296, 10987178906592, \\ \nonumber &737967598470144, 57713234231210688, 5184724381875974016,... \end{align}
Is there a Moonshine connecting these numbers with representation theory?

The multiplier systems (see Definition 2.1 below) of $\bG$ at weight $w$ are parametrised by the representations
of $\overline{\bG}$. In dimension $d$, the  space of $\overline{\bG}$-representations,
--- more precisely, the algebraic quotient of all group homomorphisms
$\overline{\bG}\rightarrow GL_d(\bbC)$ by the conjugate action of $GL_d(\bbC)$,
form a variety, and the completely reducible  representations form an open subvariety.
The connected components of this open subvariety correspond to  ordered
5-tuples $(\alpha_i;\beta_j)$ \cite{bib:West},
where $\alpha_i$ is the eigenvalue multiplicity of  $(-1)^i$ for $S$, and $\beta_j$ is that of $\xi_3^{j}$
for $U$. Of course $\alpha_0+\alpha_1=\beta_0+\beta_1+\beta_2 =d$. When $d>1$, the 
$(\alpha_i;\beta_j)$ component
is nonempty iff 
\begin{equation}\mathrm{max}\{\beta_j\}\le\mathrm{min}\{\alpha_i\}\,,\label{eq:inequal}\end{equation} 
in which case its dimension 
is $d^2+1-\sum_i\alpha_i^2-\sum_j\beta_j^2$ \cite{bib:West}. The irreducible $\overline{\bG}$-representations of dimension $d<6$ are explicitly described in
\cite{bib:TW}.
For irreducible representations, \eqref{eq:inequal} is obtained using quiver-theoretic means in \cite{bib:West}; later we recover  \eqref{eq:inequal}  within our theory.

\subsection{Definitions}

 Write $1_d$ for the $d\times d$ identity matrix, and $e_j=(0,\ldots,1,\ldots,0)^t$ its $j$th
 column.

\medskip{}
\noindent\textbf{Definition 2.1(a)} \textit{An} admissible multiplier system $(\rho,w)$
\textit{consists of some $w\in\bbC$ called the} weight and \textit{a map} $\map{\rho}{\bG}{\mathrm{GL}_{d}(\mathbb{C})}$ 
\textit{called the} multiplier, \textit{for some
positive integer $d$ called the} rank, \textit{such that:}

\smallskip\noindent\textbf{(i)} \textit{the associated automorphy factor
\[
\widetilde{\rho}_{w}\!\left(\gamma,\tau\right)=\rho\!\left(\gamma\right)
%\left(\dif{\!\left(\gamma\tau\right)}{\tau}\right)^{-w/2}
\,(c\tau+d)^w\]
satisfies, for all} $\gamma_{1},\gamma_{2}\in\bG$, \begin{equation}
\widetilde{\rho}_{w}(\gamma_{1}\gamma_{2},\tau)=\widetilde{\rho}_{w}(\gamma_{1},\gamma_{2}\tau)\,\,
\widetilde{\rho}_{w}(\gamma_{2},\tau)\,;\label{eq:mulsys}\end{equation}

\smallskip\noindent\textbf{(ii)} \textit{$\rho(1_2)$ and 
$e^{-\pi\Ii w}\rho(-1_2)$ both equal the identity matrix.}
\medskip

The conditions on $\rho(\pm 1_2)$ in (ii) are necessary
(and sufficient, as we'll see) for the existence of 
nontrivial vvmf at weight $w$.  In practise
modular forms are most important for their Fourier expansions. This is often true for vvmf,
and this is why we take $\rho$ to be  matrix-valued in Definition 2.1. Note that the
multiplier $\rho(\gamma)$ need not be unitary, and the weight $w$ need not be real.

When $w\in\bbZ$, $(\rho,w)$ is admissible iff $\rho$ is a representation of $\Gamma$
satisfying $\rho(-1_2)=e^{\pi\Ii w}$. When $w\not\in\bbZ$, $\rho$ is only a \textit{projective}
representation of $\Gamma$, and is most elegantly described in terms of the braid group
$B_3$. More precisely, $B_3=\langle\sigma_1,\sigma_2\,|\,\sigma_1\sigma_2\sigma_1=\sigma_2\sigma_1\sigma_2\rangle$ is a central extension by $\bbZ$
of $\bG$, where the surjection $B_3\rightarrow\bG$ sends 
$\sigma_1\mapsto \left({1\atop 0}{1\atop 1}\right)$ and $\sigma_2\mapsto\left({1\atop -1}{0\atop 1}
\right)$. The kernel is $\langle(\sigma_1\sigma_2\sigma_1)^4\rangle$, half of the centre of
$B_3$. Then $(\rho,w)$ is admissible iff there is a representation $\widehat{\rho}$ of
$B_3$ (necessarily unique) satisfying
\begin{equation}\widehat{\rho}(\sigma_1)= T\,,\ \ \widehat{\rho}((\sigma_1\sigma_2\sigma_1)^{-1})
= S\,,\ \ \widehat{\rho}( (\sigma_1\sigma_2)^{-2})=U\label{braid}\end{equation} 
and also $\widehat{\rho}((\sigma_1\sigma_2\sigma_1)^2)=e^{\pi\Ii w}$. Alternatively, we
see in  Lemma 3.1 below that for any  
 $w\in\bbC$, there is an admissible system $(\varmu_{w},w)$ of rank 1; then $(\rho,w)$
 is admissible iff $\overline{\varmu_w}\otimes\rho$ is a $\overline{\bG}$-representation. 
 From any of these descriptions, we see that a multiplier
$\rho$ determines the corresponding weight $w$ modulo 2.

\medskip{}
\noindent\textbf{Definition 2.2.} \textit{Let $(\rho,w)$ be an admissible multiplier system 
of rank $d$. A map $\map{\FA}{\uhp}{\mathbb{C}^{d}}$
is called a}  vector-valued modular form (vvmf) \textit{provided \begin{equation}
\FA\left(\gamma\tau\right)=\widetilde{\rho}_{w}\!\left(\gamma,\tau\right)\FA\left(\tau\right)\label{eq:modtrans}\end{equation}
 for all $\gamma\!\in\!\bG$ and $\tau\!\in\!\uhp$, and each component 
$\FA\,{}_i\left(\tau\right)$ is meromorphic throughout $\uhp^*$.
 We write $\cM^!_w(\rho)$ for the space of all} weakly holomorphic \textit{vvmf, i.e.
those holomorphic throughout $\uhp$.
%, and $\cH_w(\rho)$ for the ones holomorphic throughout $\uhp^*$.
} \medskip

Meromorphicity at the cusps is defined as usual, e.g. by a growth condition or through the
$q$-expansion given shortly.

Generic $\rho$ will have $T$ diagonalisable, which we can then insist is diagonal without
loss of generality. For simplicity, we assume throughout this paper that $T$ is
diagonal. The following theory generalises to $T$ a direct sum of Jordan blocks
(the so-called `logarithmic' case) without difficulty, other than notational awkwardness
\cite{bib:BG3,bib:Ga}. 

Assume then that $T$ is diagonal. By an \textit{exponent} $\lambda$ for $\rho$, we mean
any diagonal matrix such that $e^{2\pi\Ii\lambda}=T$, i.e. $T_{jj}=e^{2\pi\Ii\lambda_{jj}}$ for.
all $j$. An exponent is uniquely defined modulo 1. It is typical in the literature to fix the real
part of $\lambda$ to be between 0 and 1. But we learn in Theorem 3.2 below that 
often there will be better exponents to choose.

For any vvmf $\bbX$ and any exponent $\lambda$, $q^{-\lambda}\bbX(\tau)$ will be
invariant under $\tau\mapsto\tau+1$, where we write $q^\lambda=\mathrm{diag}(e^{2\pi\Ii\tau\lambda_{11}},\ldots,
e^{2\pi\Ii\tau\lambda_{dd}})$. This gives us a Fourier expansion 
\begin{equation}\label{eq:fourier}\bbX(\tau)=
q^\lambda\sum_{n=-\infty}^\infty\bbX_{(n)}\,q^n\,,\end{equation} where the coefficients $\bbX_{(n)}$ lie
in $\bbC^d$. $\bbX(\tau)$ is meromorphic at the cusp $\Ii\infty$ iff only finitely many coefficients
$\bbX_{(n)}$, for $n<0$, are nonzero.

\subsection{Local expansions}

We are now ready to describe the local expansions about any of the 3 special points $\Ii\infty,
\Ii,\xi_6$ (indeed, the same method works for any point in $\uhp^*$).

\medskip
\noindent\textbf{Lemma 2.1.} \textit{Let $(\rho,w)$ be admissible
 and $T$ diagonal. Then any $\Fa\in\cM^!_{w}({\rho})$ obeys}
\begin{equation}\label{expansions}
{\Fa}(\tau)=q^{\lambda}\sum_{n=0}^\infty\Fa_{(n)}\,q^n%\in q^{\lambda}\,V_T\otimes_\bbC\bbC[[q]]
%\,.\label{eq:csp}\\
%Moreover, there are invertible matrices $P_{2},P_{3}$ % and diagonal integer matrices $\lambda_{(\I)},\lambda_{(\zeta)}$ 
%such that, when we write
%is an invertible matrix $P_{(z)}$, a partition $\sum_i \ell_{(z)i}=d$ of $d$, and 
%a diagonal matrix 
%$\lambda_{(z)}=\mathrm{diag}(\lambda_11_{\ell_{(z)1}},\ldots,\lambda_s1_{\ell_{(z)s}})$,
%such that $P_{(z)}\rho(\gamma_{(z)})P_{(z)}^{-1}=e^{2\pi\I(\bN(\vec{\ell}_{(z)})+\lambda_{(z)})}$ and, writing 
%$\widetilde{\Fa}_{(z)}(\tau_{(z)})=j_{(z)}(w;\tau)P_{(z)}\Fa(\tau)$ 
%such that for $i=2$ and $3$ we have
%\begin{align}
%\widetilde{\Fa}_{(z)}(e^{2\pi\I/\eps(z)}\tau_{(z)})=&\,e^{-\pi\I w/\eps(z)}
%P_{(z)}\rho(\gamma_{(z)})P_{(z)}^{-1}
%\widetilde{\Fa}_{(z)}(\tau_{(z)})\quad\mathrm{for}\ z\in\uhp\,,\label{eq:ellfpa}\\
%\widetilde{\Fa}_{(z)}(\tau_{(z)}+h(z))=&\,P_{(z)}\rho(\gamma_{(z)})
%P_{(z)}^{-1}\widetilde{\Fa}_{(z)}(\tau_{(z)})\qquad\mathrm{for}\ z\in\cC\,.
%\label{eq:cspa}
%\end{align}
%Moreover, for any $z\in\uhp^*$,
%\begin{equation}
=\mathrm{j}_{2}(w;\tau)^{-1}\sum_{n=0}^\infty\Fa_{[n]}\,\tau_{2}^n
%\in \bbC^d[[\tau_{2}]]\label{eq:ell2}\\ 
=\mathrm{j}_{3}(w;\tau)^{-1}\sum_{n=0}^\infty\Fa_{\langle n\rangle}\,\tau_{3}^n
%\in \bbC^d[[\tau_{3}]]\label{eq:ell3}
\,,\end{equation}
\textit{for some exponent $\lambda$. These converge for $0<|q|<1$ and $|\tau_i|<\epsilon_i$. Also,} 
\begin{equation}e^{\pi\Ii w/2}S\Fa_{[n]}=(-1)^{n}\Fa_{[n]}\qquad\mathrm{and}\qquad e^{2\pi\Ii w/3}
U\Fa_{\langle n\rangle}=\xi_3^{n}\Fa_{\langle n\rangle}\,.\label{SUtrans}\end{equation}

%\smallskip\noindent\textbf{(b)} \textit{
%It is possible to choose these $P_{(z)}$ so that, whenever $z,z'$ lie in the
%same $\bG$-orbit, say $z'=\gamma z$ for $\gamma=\left({a\atop c}{b\atop d}\right)\in\bG$, then 
%\begin{equation} \widetilde{\Fa}_{(z')}(\tau_{(z')})=f(z)\,
%\widetilde{\Fa}_{(z)}((\gamma^{-1}\tau)_{(z)})\,,\label{eq:ellorb}\end{equation}
%where $f(z)=(cz+d)^w$ for $z\in\uhp$, $f(z)=(cz+d)^{-w}$ for $z,z'\in\bbR$,
%and $f(z)=c^{-w}$ for $z=\infty,z'\in\bbR$.}
\medskip

The existence of the $q$-series follows from the explanation at the end of last subsection.
For $\tau=\Ii$ ($\tau=\xi_6$ is identical), j$_2(w;\tau)\,\bbX(\tau)$ is holomorphic in the disc
$|\tau_2|<\epsilon_2$ and so has a Taylor expansion. The transformation \eqref{SUtrans}
can be seen by direct calculation, but will be trivial once we know Lemma 3.1 below.

We label components by $\bbX_{(n)\, i}$ etc. A more uniform notation would have been
to define e.g. $q_2=\tau^2_2$, find a matrix $P_2$  and  an `exponent matrix' $\lambda_2$ whose
diagonal entries lie in $\frac{1}{2}\bbZ$, such that $P_2SP_2^{-1}=e^{2\pi\Ii\lambda_2}$ and
$\bbX=\mathrm{j}_2(w;\tau)^{-1}P_2^{-1}q_2^{\lambda_2}\sum_{n=0}^\infty \widetilde{\bbX}_{[n]}
q_2^n$. For most purposes
the simpler \eqref{expansions} is adequate, but see \eqref{xiexp2} below.

\subsection{Differential operators}

Differential equations have played a large role in the theory of vvmf.
The starting point is always the modular derivative 
\begin{equation}\Der_wf=\frac{1}{2\pi\Ii}\frac{\D}{\D\tau}-\frac{w}{12}E_2=q\frac{\D}{\D q}-
\frac{w}{12}E_2\,,\end{equation}
where $E_2(\tau)=1-24q-72q^2-\cdots$ is the quasi-modular Eisenstein function.
 Note that  $\Der_{12}$ kills the discriminant form $\Delta(\tau)=\eta(\tau)^{24}$. This 
$\Der_w$ maps $\cM^!_w(\rho)$ to $\cM^!_{w+2}(\rho)$. It is a derivation in the sense that
if $f\in\frakm$ is weight $k$ and $\bbX\in\cM_w^!(\rho)$, then $\Der_{k+w}(f\bbX)=
\Der_k(f)\,\bbX+f\,\Der_w(\bbX)$.
We write $\Der_w^j=
\Der_{w+2j-2}\circ\cdots\circ\Der_{w+2}\circ\Der_w$.

There are several different applications of differential equations to modular forms --- some are
reviewed in \cite{bib:Za}. But outside of our work, the most influential for the theory of vvmf (see 
e.g. \cite{bib:AM,bib:Mil,bib:Mas,bib:Mar}) has been the differential equation
coming from the Wronskian (see e.g. \cite{bib:Mas} for the straightforward proof):

\medskip\noindent\textbf{Lemma 2.2(a)} \textit{Let $(\rho,w)$ be admissible. For $\Fa\in\cM^!_w(\rho)$, define}
\begin{equation}\mathrm{Wr}(\Fa):=\mathrm{det}\left(\begin{matrix}\Fa_1&\frac{\mathrm{d}}{2\pi\Ii\mathrm{d}\tau}\Fa_1&\cdots&(\frac{\mathrm{d}}{2\pi\Ii\mathrm{d}\tau})^{d-1}\Fa_1\\
\vdots&\vdots&&\vdots\\ \Fa_d&\frac{\mathrm{d}}{2\pi\Ii\mathrm{d}\tau}\Fa_d&\cdots&(\frac{\mathrm{d}}{2\pi\Ii\mathrm{d}\tau})^{d-1}\Fa_d\end{matrix}\right)
=\mathrm{det}\left(\begin{matrix}\Fa_1&\Der_w\Fa_1&\cdots&\Der_w^{d-1}\Fa_1\\
\vdots&\vdots&&\vdots\\ \Fa_d&\Der_w\Fa_d&\cdots&\Der_w^{d-1}\Fa_d\end{matrix}\right)\,.
\nonumber\end{equation}
\textit{Then Wr$(\Fa)(\tau)\in\cM^!_{(d+w-1)d}(\mathrm{det}\,\rho)$. 
If the coefficients of $\Fa$ are linearly independent over 
$\bbC$, then the function Wr$(\Fa)(\tau)$ is nonzero.}\smallskip

\noindent\textbf{(b)} \textit{Given admissible $(\rho,w)$ and $\bbX\in\cM^!_w(\rho)$,
 define an operator $L_\Fa$ on the space of all functions $y$ meromorphic on
$\uhp^*$, by
\begin{equation}\label{eq:ODE}
L_\Fa=\mathrm{det}\left(\begin{matrix}y&\Der_wy&\cdots&\Der_w^{d}y\\\Fa_1&\Der_w\Fa_1&\cdots&\Der_w^d\Fa_1\\
\vdots&\vdots&&\vdots\\ \Fa_d&\Der_w\Fa_d&\cdots&\Der_w^{d}\Fa_d\end{matrix}\right)=
\sum_{l=0}^dh_l(\tau)\,\Der_{w}^ly\,,\end{equation}
where $h_d=\mathrm{Wr}(\Fa)$ and each $h_l$ is a (meromorphic scalar) modular form of weight $(d+w+1)d-2l$
with multiplier det$\,\rho$. Then $L_\Fa\,\Fa_i=0$ for
all components $\Fa_i$ of $\Fa$. Conversely, when the components of $\Fa$ are
linearly independent, the solution space to $L_\Fa y=0$ is  $\mathrm{Span}_\bbC \{\Fa_i\}$.} 
\medskip

In our theory, the differential equation \eqref{eq:ODE} plays a minor role. More important
are the differential operators which don't change the weight:
\begin{equation}\nabla_{1,w}=\frac{E_{4}E_6}{\Delta}\Der_w\,,\ \ \ \nabla_{2,w}=\frac{E_4^2}{\Delta}
\Der_w^2\,,\ \ \ \nabla_{3,w}=\frac{E_6}{\Delta}\Der_w^3\,.\end{equation}
Each $\nabla_{i,w}$ operates on $\cM^!_w(\rho)$, and an easy calculation shows that 
any $f\Der_w^j$ for $f\in\cM^!_{-2j}(1)$
is a polynomial in these three $\nabla_{i,w}$ with coefficients in $\bbC[J]$. Conversely, $\nabla_{3,w}$
is not in $\bbC[J,\nabla_{1,w},\nabla_{2,w}]$ and $\nabla_{2,w}$ is not in $\bbC[J,\nabla_{1,w}]$.
The reason $\nabla_{2,w}$ and $\nabla_{3,w}$ are needed is because of the elliptic
points of order 2 and 3 --- this is made explicit in the proof of Proposition 3.2. It is crucial to our theory that $\cM^!_w(\rho)$ is a module
over $\bbC[J,\nabla_{1,w},\nabla_{2,w},\nabla_{3,w}]$. 

We sometimes drop the subscript
$w$ on $\Der_w$ and $\nabla_{i,w}$ for readability.

\section{Our main results}

\subsection{Existence of vvmf}

The main result (Theorem 3.1) of this subsection is the existence proof of vvmf for any
$(\rho,w)$. 
As a warm-up, let us show there is a (scalar) modular form of every complex weight, and compute its multiplier.

\medskip\noindent\textbf{Lemma 3.1.} \textit{For any $w\in\bbC$, there is a weakly holomorphic  modular
form $\Delta^{w}(\tau)=q^w(1-24w q+\cdots)$ of  weight $12w$, nonvanishing everywhere except at the cusps. The multiplier $\varmu_{12w}$, in terms of the braid group $B_3$  (recall \eqref{braid}), is}
\begin{equation}\widehat{\varmu}_{12w}(\sigma_1)=\widehat{\varmu}_{12w}(\sigma_2)=\exp(2\pi\Ii w)\,.
%,\quad\varmu_{12w}((\sigma_1\sigma_2\sigma_1)^{-1})=\exp(-6\pi\Ii w)\,,\quad
%\varmu_{12w}((\sigma_1\sigma_2)^{-2})=\exp(-8\pi\Ii w)\,.
\label{nuw}\end{equation}

\smallskip\noindent\textit{Proof.} First note that the discriminant form $\Delta(\tau)$
is holomorphic and nonzero in the simply-connected domain  $\bbH$, and so has
a logarithmic derivative Log$\,\Delta(\tau)$ there. Hence $\Delta^{w}(\tau):=\exp(w\mathrm{Log}\,
\Delta(\tau))$ is well-defined and holomorphic throughout $\bbH$. It is easy to verify that $\Delta^{w}(\tau)$ satisfies the differential equation
\begin{equation} 
\frac{1}{2\pi\Ii}\frac{\D f}{\D\tau}=wE_2(\tau)\,f(\tau)\label{Deltadiffeq}\end{equation}
--- indeed, this simply reduces to the statement that $E_2(\tau)$ is the logarithmic derivative of $\Delta(\tau)$.
Therefore \textit{any}  solution to \eqref{Deltadiffeq} is a scalar multiple of $\Delta^{w}(\tau)$.

Now, fix $\gamma\in\Gamma$. Then $f(\tau)= (c\tau+d)^{-12w}\Delta^{w}(\gamma\tau)$ exists
and is holomorphic throughout $\bbH$ for the same reason. Note that $f(\tau)$ also satisfies
the differential equation \eqref{Deltadiffeq}:
\begin{align} \frac{1}{2\pi\Ii}\frac{\D}{\D\tau}f(\tau)=&\frac{-6cw}{\pi\Ii}(c\tau+d)^{-12w}(c\tau+d)^{-1}
\Delta^w(\gamma\tau)\nonumber\\&+(c\tau+d)^{-12w}(c\tau+d)^{-2}wE_2(\gamma\tau)\Delta^w(\gamma\tau)=wE_2(\tau)\,f(\tau)\,,\nonumber\end{align}
using quasi-modularity of $E_2(\tau)$, and thus $f(\tau)=\varmu \Delta^w(\tau)$ throughout $\bbH$, for some constant $\varmu=\varmu(\gamma)\in\bbC$.

The final step needed to verify that $\Delta^w(\tau)$ is a modular form of weight $12w$ is that it behaves well at the
cusps. By the previous paragraph it suffices to consider $\Ii\infty$. But there $\Delta^w(\tau)$ has
the expansion $\Delta^{w}(\tau)=q^w(1-24w q+\cdots)$, up to a constant factor which we can take to be 1.
As all cusps lie in the $\Gamma$-orbit of $\Ii\infty$, we see that $\Delta^w(\tau)$ is indeed holomorphic at all cusps.

In general, $\varmu_{12w}$ can be interpreted as a representation $\widehat{\varmu}$ of $B_3$.
The expansion tells us that $\Delta^w(\tau+1)=\exp(2\pi\Ii w)\Delta^w(\tau)$, so we have
$\widehat{\varmu}_{12w}(\sigma_1)=\exp(2\pi\Ii w)$. From the familiar $B_3$ presentation we see that any 
one-dimensional representation of $B_3$ takes the same value on both generators $\sigma_i$,
so $\widehat{\varmu}_{12w}$ is determined. \qed% \qquad\textbf{QED} to Lemma 3.1
\medskip 

For admissible $(\rho,w)$, Lemma 3.1 says the matrices $e^{\pi\Ii w/2} S$ and $e^{2\pi\Ii w/3}U$ 
have order 2 and 3 respectively. Write $\alpha_j{(\rho,w)}$ for the multiplicity of $(-1)^j$ as
an eigenvalue of  $e^{\pi\Ii w/2} S$, and $\beta_j{(\rho,w)}$ for the  eigenvalue multiplicity of
 $\xi_3^j$ for $e^{2\pi\Ii w/3}U$.

\medskip\noindent\textbf{Theorem 3.1.} \textit{Let $(\rho,w)$ be admissible and $T$ diagonal.
 Then there is a $d\times d$ matrix $\Phi(\tau)$ and exponent $\lambda$  such that
the columns of $\Phi$ lie in $\cM^!_{w}(\rho)$, and}
\begin{equation}\label{eq:phiexp}
\Phi(\tau)= q^\lambda\,\bF(q)=q^\lambda\sum_{n=0}^\infty F_{(n)}q^n\,,
\end{equation}
\textit{where the matrix $\bF(q)$ is holomorphic and invertible
 in a neighbourhood of $q=0$.}\medskip

\noindent\textit{Proof.} Consider any representation $\rho'$ of $\pi_1(\bbP^1\setminus\{0,1,\infty\})
\cong F_2$, the free group with 2 generators.   R\"orhl's solution   \cite{bib:Roh} to the \textit{Riemann--Hilbert problem} \cite{bib:Bolicm,bib:Bolrev}
(see also \cite{bib:Kim}) says that there exists a Fuchsian differential equation 
\begin{equation}\label{eq:fuchde}
\frac{\D}{\D z}{\Psi}(z)={\Psi}\left(\frac{A_1}{z}+\frac{A_2}{z-1}+\frac{B}{z-b}\right)
\end{equation}
on the Riemann sphere $\bbP^1$ whose monodromy is given by $\rho'$ --- i.e. the monodromy
corresponding to a small 
circle about  0 and 1, respectively, equals the value of $\rho'$ at the corresponding loops in
the fundamental group.  There will also be a simple pole at $\infty$, with residue $A_\infty:=-A_0-A_1-B$.
The $B$ term in \eqref{eq:fuchde} corresponds to an apparent singularity; it can be dropped
if the monodromies about 0 or 1 or $\infty$ have finite order \cite{bib:Bolicm,bib:Bolrev} (which  will happen for the $\rho'$ of interest to us).
About any of these 3 or 4 singular points $c\in\{0,1,\infty,b\}$,  Levelt \cite{bib:Lev} proved that a solution ${\Psi}
(z)$ to such a differential equation has the form 
\begin{equation}\Psi(z)=P_c^{-1}\widetilde{z}^{N_c}\widetilde{z}^{\lambda_c}\mathrm{F}_c(z)
\label{levelt}\end{equation} where ${N}_c$
is nilpotent, $\lambda_c$ is diagonal, $\mathrm{F}_c(z)$ is holomorphic and holomorphically 
invertible about $z=c$, and $\widetilde{z}$ is a coordinate on the universal cover of a small
disc punctured at $c$. By \cite{bib:Kim}, we may take ${N}_c+\lambda_c$ to be conjugate to $A_c$.

Now suppose we have a representation of the free product $\overline{\bG}\cong 
\bbZ_2*\bbZ_3$. This  is a homomorphic image of the free group $F_2$, so we can lift $\rho'$
to $\pi_1(\bbP^1\setminus\{0,1,\infty\})$. For us, $\bbP^1$ is $\overline{\bG}\backslash\uhp^*$, with (smooth) global coordinate $z=J(\tau)/1728$.
The monodromy at $z=0$ and $z=1$ has finite order 2 and 3 respectively, so we won't need
the apparent singularity $b$. The point $\infty$ corresponds to $\Ii\infty$ (or rather its $\Gamma$-orbit),
where $q_\infty=q$, and we have $N_\infty=0$ and $P_\infty=1_d$ since $T$ is diagonal; in this case Levelt's equation \eqref{levelt} 
reduces to \eqref{eq:phiexp}. The singularities 0 and 1 correspond to the order 2 and 3 
elliptic points $\Ii$ and $\xi_6$; for them, $N=0$, $z=\tau_i$ locally looks like $(\tau-\Ii)^2$ and 
$(\tau-\xi_6)^3$, and the diagonal elements of $\lambda$ lie in $\frac{1}{2}\bbZ$ and 
$\frac{1}{3}\bbZ$. Then Levelt's equation \eqref{levelt} says $\Psi(J(\tau)/1728)$ is meromorphic at $\Ii$ and $\xi_6$
(hence their $\Gamma$-orbits). It is automatically holomorphic at all other points.
   
 The desired matrix is $\Phi(\tau)={\Psi}(J(\tau))\,\Delta^{w/12}(\tau)J(\tau)^m(J(\tau)-1)^n$,
 where  $m,n\in\bbZ_{\ge 0}$ are taken large enough to kill any poles at $z=0$ and 1 (i.e. at 
 all elliptic points), and $\Psi$ corresponds to the $\overline{\Gamma}$-representation
 $\rho'=\varmu_{-w}\otimes \rho$. 
 \qed%\textbf{QED} to Theorem 3.1
 \medskip

The proof generalises without change to nondiagonalisable $T$, and to any other genus-0 Fuchsian group of the first kind \cite{bib:Ga}.
Theorem 3.1 is vastly more general than previous vvmf existence proofs.
Previously (see \cite{bib:KM,bib:KMlog}), 
existence of vvmf was only established for $\bG$ with real weight $w$, and
requiring in addition that the eigenvalues of $T$ all
have modulus 1. Their proof used Poincar\'e series; the difficult step there
is to establish convergence, and that is what has prevented their methods
to be generalised. That analytic complexity was handled here by R\"ohrl's argument.

\subsection{Mittag-Leffler}

Let $(\rho,w)$ be admissible and $T$ diagonal. In this subsection we study the principal part
map and calculate its index. This is fundamental to our theory. As always in this paper, the generalisation to
nondiagonalisable $T$ and to other genus-0 Fuchsian groups is straightforward \cite{bib:BG3}. 

Given any exponent $\lambda$ and any $\bbX(\tau)\in\cM^!_w(\rho)$, we have the
$q$-expansion \eqref{eq:fourier}.
 Define the  \textit{principal part map} $\cP_\lambda:\cM^!_w(\rho)\rightarrow \bbC^d[q^{-1}]$ by
\begin{equation} \cP_\lambda(\bbX)=\sum_{n\le 0}\bbX_{(n)}\,q^n\,.\end{equation}
When we want to emphasise the domain, we'll write this $\cP_{\lambda;(\rho,w)}$.

\medskip
\noindent\textbf{Theorem 3.2.} \cite{bib:BG3} \textit{Assume $(\rho,w)$ is admissible, and $T$ is
diagonal. Recall the eigenvalue multiplicities $\alpha_j=\alpha_j{(\rho,w)}$ and $\beta_j=\beta_j{(\rho,w)}$
from Section 3.1.}

\smallskip\noindent\textbf{(a)} \textit{For any exponent $\lambda$, $\cP_\lambda:\cM^!_w(\rho)\rightarrow
\bbC^d[q^{-1}]$ has finite-dimensional kernel and cokernel, and the index is}
\begin{equation}\label{index} \mathrm{dim\ ker}\,\cP_{\lambda}-\mathrm{dim\ coker}\,\cP_\lambda
=-\mathrm{Tr}\,\lambda+c_{(\rho,w)}\,,\end{equation}
\textit{for}\begin{align}
c_{(\rho,w)}=&\,\frac{(w-7)d}{12}+\frac{e^{\pi\Ii w/2}}{4}\mathrm{Tr}\,{S}+
\frac{2}{3\sqrt{3}}\mathrm{Re}\left(e^{\frac{-\pi\Ii}{6}-\frac{2\pi\Ii w}{3}}
\mathrm{Tr}\,{U}\right)\:\nonumber\\
=&\,\frac{wd}{12}-\frac{\alpha_1}{2}-\frac{\beta_1+2\beta_2}{3}\,.\label{cvalue}
\end{align}

\smallskip\noindent\textbf{(b)} \textit{There exist exponents $\lambda$ for which $\cP_\lambda:
\cM^!_w(\rho)\rightarrow\bbC^d[q^{-1}]$ is a vector space isomorphism.}\medskip

By a \textit{bijective exponent} we mean any exponent $\lambda$ for which $\cP_\lambda:\cM^!_w(\rho)\rightarrow
\bbC^d[q^{-1}]$ is an isomorphism. Of course by \eqref{index} its trace $\sum_j\lambda_{jj}$ must equal
$c_{(\rho,w)}$, but the converse is not true as we will see.

For example, for the trivial 1-dimensional representation, $T=1$ so an exponent is just an
integer. Here, $ \cM^!_0(1)=\bbC[J]$ and $c_{(1,0)}=0$. The map $\cP_1$ is injective but not
surjective (nothing has principal part 1), while  $\cP_{-1}$ is surjective but not injective (it kills all constants). For another example, taking $\rho=\varmu_{-2}$ (the multiplier 
of $\eta^{-4}$), we have $c_{(\varmu_{-2},0)}=-7/6$.

It is standard in the literature to restrict from the start to exponents satisfying $0\le \lambda_{jj}
<1$. However, such $\lambda$
are seldom bijective. It is rarely wise  to casually throw away a freedom.

Theorem 3.2(b) first appeared in \cite{bib:BG2}, though for restricted $(\rho,w)$, and with the
erroneous claim that $\lambda$ is bijective iff Tr$\,\lambda=c_{(\rho,w)}$.
The deeper part of Theorem 3.2 is  the index formula, which is new. We  interpret it later as Riemann--Roch,
and obtain from it dimensions of spaces of holomorphic vvmf.

The right-side of \eqref{cvalue} is always integral. To see that, take $w=0$ and note
\begin{equation}\label{trmod6}
\exp(2\pi\Ii\,\mathrm{Tr}\,\lambda)=\mathrm{det}\,T=\mathrm{det}\,S\,\,\mathrm{det}\,U^{-1}=
(-1)^{\alpha_1}{\xi}_3^{-\beta_1+\beta_2}\,.\end{equation}

Fix an admissible $(\rho,w)$ with diagonal $T$, and a  bijective exponent $\Lambda$. As a vector 
space over $\bbC$, $\cM^!_w(\rho)$ has a basis 
\begin{equation}\label{eq:basis}\bbX^{(j;n)}(\tau)
=\cP_\Lambda^{-1}(q^{-n}e_j)=q^\Lambda\left(q^{-n}e_j+
\sum_{m=1}^\infty\bbX^{(j;n)}_{(m)}q^m\right)\,,\end{equation}
 where $e_j=(0,\ldots,1,\ldots,0)^t$ and $n\in\bbZ_{\ge 0}$. We describe
next subsection an effective way to find all these $\bbX^{(j;n)}(\tau)$, given the $d^2$ coefficients \begin{equation}\label{eq:char}
\chi_{ij}:=\bbX^{(j;0)}_{(1)\, i}\in\bbC\,.\end{equation}
 We learn there that $\cM^!_w(\rho)$ is under total control once a  bijective exponent
$\Lambda$ and its corresponding matrix $\chi=\chi(\Lambda)$ are found.

Recall that for fixed $\rho$, the weight $w$ 
is only determined mod 2, i.e. $(\rho,w)$ is admissible iff $(\rho,w+2k)$ is, for any $k\in\bbZ$. 
We find from the definition of the $\alpha_i$ and $\beta_j$ that
\begin{equation}\label{alphabet}
 \alpha_j{(\rho,w+2k)}=\alpha_{j+k}{(\rho,w)}\,,\qquad\beta_j{(\rho,w+2k)}=\beta_{j+k}(\rho,w)\,.
 \end{equation}  
 Plugging this into \eqref{cvalue}, we obtain the trace of a bijective exponent for $(\rho,w+2k)$:
 \begin{equation}\label{tracek}
 c_{(\rho,w+2k+12l)}=c_{(\rho,w)}+ld+\left\{\begin{matrix}0&&&&\mathrm{if}\ k=0\\ \alpha_1-\beta_0&&&&
 \mathrm{if}\ k=1\\ \beta_2&&&& \mathrm{if}\ k=2\\ \alpha_1&&&& \mathrm{if}\ k=3\\ \beta_1+\beta_2&&&&
 \mathrm{if}\ k=4\\ \alpha_1+\beta_2&&&& \mathrm{if}\ k=5\end{matrix}\right.\ \ .\end{equation}
 
 \medskip\noindent\textit{Sketch of proof of Theorem 3.2}. (see \cite{bib:BG3} for the complete proof and its generalisation). First of all, it is easy to show that if the
real part of an exponent $\lambda$ is sufficiently large, then $\cP_\lambda$ is necessarily
injective. This implies that the kernel of any $\cP_{\lambda'}$ is finite-dimensional, in
fact $\mathrm{dim\ ker}\ \cP_{\lambda'}\le \sum_j\mathrm{max}\{\lambda_{jj}-\lambda'_{jj},0\}$
for any $\lambda$ with $\cP_\lambda$ injective. 

Let $M$ be the 
$\bbC[J]$-span of the columns of the matrix $\Phi$ of Theorem 3.1. Let $\lambda_M$ be
the corresponding exponent appearing in \eqref{eq:phiexp}. Then invertibility of F$(q)$ implies
invertibility of $F_{(0)}$, which in turn implies $\cP_{\lambda_M}:M\rightarrow\bbC^d[q^{-1}]$ is surjective, and hence so is $\cP_{\lambda_M}:\cM^!_w(\rho)
\rightarrow\bbC^d[q^{-1}]$. This implies any cokernel is also finite-dimensional.

A bijective exponent can be obtained by starting from $\lambda_M$, and recursively increasing
one of its entries by $+1$ so that something in the kernel of the previous $\cP$ is no longer
in the kernel of the new $\cP$. This implies the range is unchanged. This process must terminate, by finiteness of the kernel. 

The index formula \eqref{index}, for some value of $c_{(\rho,w)}$, is now computed by showing that, whenever $\lambda'\ge\lambda$,
the index of $\cP_\lambda$ minus that of $\cP_{\lambda'}$ equals Tr$\,(\lambda'-\lambda)$.
The constant $c_{(\rho,w)}$ is computed next subsection. \qed%\qquad\textbf{QED} to Thm. 3.2
\medskip

\noindent\textbf{Corollary 3.1.} \textit{Suppose $(\rho,w)$ is admissible, $T$ is diagonal,
 and $\rho$ is completely reducible with  no 1-dimensional subrepresentation. 
Then setting $\epsilon=0,1$ for $d$ even, odd respectively, we obtain the bounds} 
\begin{equation}
\frac{wd}{12}-d+\frac{\epsilon}{4}\le c_{(\rho,w)}\le\frac{wd}{12}-\frac{5d}{12}
-\frac{\epsilon}{4}\,.\label{eq:Trbound}
\end{equation} 

\noindent\textit{Proof.}
Assume without loss of generality that weight $w=0$. Then \eqref{cvalue} tells us that $c_{(\rho,w)}
=-\alpha_1/2-(\beta_1+2\beta_2)/3$, where $\alpha_i=\alpha_i(\rho,w),
\beta_j=\beta_j(\rho,w)$.

First, let's try to maximise $c_{(\rho,w)}$, subject to the inequalities \eqref{eq:inequal}.
Clearly, $c_{(\rho,w)}$ is largest when $\alpha_0\ge\alpha_1\ge
\beta_0\ge\beta_1\ge\beta_2$.  In fact, we should take $\beta_0$ as large as possible, i.e.
 $\alpha_1=\beta_0$. Then our formula for $c_{(\rho,w)}$ simplifies to $-\frac{d}{2}+\frac{
 \beta_1-\beta_2}{6}$. It is now clear this is maximised by $\beta_2=\epsilon$ and $\beta_0
 =\beta_1=\frac{d-\epsilon}{2}$, which recovers the upper bound in \eqref{eq:Trbound}.
  Similarly, the lower
bound in \eqref{eq:Trbound} is realised by $\alpha_0=\beta_2=\beta_1=(d-\epsilon)/2$. 
\qed%\textbf{QED}
\medskip

The question of which exponents $\lambda$ with the correct trace are bijective,  can be 
subtle, though we see next that for
generic $\rho$ the trace condition Tr$\,\lambda=c_{(\rho,w)}$ is also 
sufficient. With this in mind,  define the $\infty\times\infty$ complex matrix $\mathcal{X}=(\mathcal{X}^{(j;n)}_{(m)\, i})$
built from the $m$th coefficient of the $i$th component of the basis \eqref{eq:basis}. 
For any $\ell\in\bbZ^d$ define a $\left(\sum_i\mathrm{max}\{\ell_i,0\}\right)\times\left(\sum_j
\mathrm{max}\{-\ell_j,0\}\right)$ submatrix $\mathcal{X}(\ell)$ of $\mathcal{X}$ by
restricting to the rows $(i;m)$ with $0<m\le\ell_i$ and columns $(j;n)$ with
$0\le n<-\ell_j$.

\medskip\noindent\textbf{Proposition 3.1.} \textit{Let $(\rho,w)$ be 
admissible, $T$ diagonal, and $\Lambda$ bijective. Then an exponent $\lambda$ is also
bijective iff the matrix $\mathcal{X}(\lambda-\Lambda)$ is invertible.} \medskip 

The effectiveness of this test will be clear next subsection, where we explain how
to compute the $\bbX^{(j;n)}$ and hence the submatrices $\mathcal{X}(\ell)$. Of course, 
invertibility forces $\mathcal{X}(\lambda-\Lambda)$ to be square, i.e.
$\sum_i\ell_i=0$, i.e. that Tr$\,\lambda=\mathrm{Tr}\,\Lambda$.

To prove Proposition 3.1, observe that the spaces $\mathrm{ker}\,\FE_{\lambda}$ and $\mathrm{null}\,\mathcal{X}
(\ell)$ are isomorphic, with $v\in\mathrm{null}\,\mathcal{X}(\ell)$ identified
with $\sum_{j}\sum_{n=0}^{-\ell_j-1}v_{(j;n)}\Fa^{(j;n)}$
(the nullspace null$\,M$ of a $m\times n$ matrix $M$
is all $v\in\bbC^n$ such that $Mv=0$). Similarly, 
the spaces $\mathrm{coker}\,\FE_{\lambda}$ and $\mathrm{null}\,\mathcal{X}(\ell)^t$ are isomorphic.

For example, given a bijective exponent $\FB$, Proposition 3.1 
says that $\chi_{ij}\ne  0$ iff  $\FB+e_i-e_j$ is bijective, where $\chi$ is defined in \eqref{eq:char}.

\subsection{Birkhoff--Grothendieck and Fuchsian differential equations}
   
As mentioned in the introduction, one of our   (vague) heuristics is to think of vvmf as
meromorphic sections of holomorphic bundles over $\bbP^1$. But the
Birkhoff--Grothendieck Theorem says that such a bundle is a direct sum of line bundles.
If taken literally, it would say that, up to a change $P$ of basis, each component $\bbX_i(\tau)$
of any $\bbX\in\cM^!_w(\rho)$ would be a \textit{scalar} modular form of weight $w$ (and some
multiplier) for $\Gamma$. This is absurd, as it only happens when $\rho$ is equivalent to
a sum of 1-dimensional projective representations. The reason Birkhoff--Grothendieck cannot
be applied here is that we do \textit{not} have a bundle (in the usual sense) over $\bbP^1$ --- 
indeed, our space $\Gamma\backslash\bbH^*$ has three singularities.

Nevertheless, Theorem 3.3(a) below says Birkhoff--Grothendieck still holds in spirit.

\medskip\noindent\textbf{Theorem 3.3.} \textit{Let $(\rho,w)$ be admissible, $T$
diagonal, and $\Lambda$ bijective.} 
   
\smallskip\noindent\textbf{(a)} \textit{$\cM^!_w(\rho)$ is a free $\bbC[J]$-module of rank $d=
\mathrm{rank}(\rho)$. Free generators are $\bbX^{(j;0)}(\tau)$ (see \eqref{eq:basis}).}

\smallskip\noindent\textbf{(b)} \textit{Let ${\Xi}(\tau)=q^\Lambda(1_d+\chi q+\sum_{n=2}^\infty{\Xi}_{(n)}\,q^n)$ 
be the $d\times d$ matrix whose columns  are the $\bbX^{(j;0)}(\tau)$. Then}
\begin{equation}
\frac{E_4(\tau)\,E_6(\tau)}{\Delta(\tau)}\Der_{w}{\Xi}(\tau)={\Xi}(\tau)\,\{(J(\tau)-984)\Lambda_w+\chi_w+[
\Lambda_w,\chi_w]\}\,,\label{fundeq1}\end{equation}
\textit{where $\Lambda_w:=\Lambda-\frac{w}{12}1_d$, $\chi_w:=\chi+2w1_d$, and $[\cdot,\cdot]$ denotes the usual bracket.}  
 
\smallskip\noindent\textbf{(c)} \textit{Assume weight $w=0$. The multi-valued function $\widetilde{{\Xi}}(z):={\Xi}(\tau(z))$,
where $z(\tau)=J(\tau)/1728$, obeys the Fuchsian differential equation}
\begin{eqnarray}&\frac{\mathrm{d}}{\mathrm{d}z}\,\widetilde{{\Xi}}(z)=\widetilde{{\Xi}}(z)\left(\frac{\cA_2}{z-1}+\frac{\cA_3}{z}\right)\,,\label{fundeq2}\\
\cA_2=\frac{31}{72}&\Lambda+\frac{1}{1728}(\chi+[\Lambda,\chi])\,,\ \
\cA_3=\frac{41}{72}\Lambda-\frac{1}{1728}(\chi+[\Lambda,\chi])\label{A2A3}\end{eqnarray}
\textit{(recall \eqref{eq:fuchde}). Moreover, $\cA_2,\cA_3$ are diagonalisable, with eigenvalues in $\{0,\frac{1}{2}\}$ and
$\{0,\frac{1}{3},\frac{2}{3}\}$, respectively.} 
   
\medskip\noindent\textit{Sketch of proof} (see \cite{bib:BG3} for the complete proof and
generalisation). The basis vvmf  $\bbX^{(j;n)}(\tau)$
exist by surjectivity of $\cP_\Lambda$. To show $\cM_w^!(\rho)$ is generated over $\bbC[J]$ by the $\bbX^{(i;0)}$, follows from an elementary induction on $n$:  if the $\bbX^{(i;m)}(\tau)$
all lie in $\sum_l\bbC[J]\bbX^{(l;0)}(\tau)$ for all $i$ and all $m<n$, then $\bbX^{(j;n)}(\tau)\in
J(\tau)\bbX^{(j;n-1)}(\tau)+\sum_l\bbC[J]\bbX^{(l;0)}(\tau)\subseteq \sum_l\bbC[J]\bbX^{(l;0)}(\tau)$, using the
fact that $\cP_\Lambda$ is injective. That
these generators are free, follows by noting the determinant of ${\Xi}(\tau)$ has a nontrivial
leading term (namely $q^\lambda$) and so is nonzero.

 The
columns of $\nabla_{1,w}{\Xi}(\tau)$ also lie in $\cM^!_w(\rho)$, and so $\nabla_{1,w}{\Xi}(\tau)={\Xi}(\tau) \,D(J)$
for some $d\times d$ matrix-valued polynomial $D(J)$, which can be determined
by Theorem 3.2 by comparing principal parts. \eqref{fundeq2} follows directly from
\eqref{fundeq1}, by changing variables. Because $e^{2\pi\Ii \cA_2}$ and $e^{2\pi\Ii \cA_3}$ must be conjugate to $S$ and $U$,
respectively, then $\cA_2$ and $\cA_3$ are diagonalisable with eigenvalues in $\frac{1}{2}\bbZ$
and $\frac{1}{3}\bbZ$. But \eqref{levelt} implies that none of these eigenvalues can be negative
(otherwise holomorphicity at $\tau=\Ii$ or $\tau=\xi_6$ would be lost). None of these eigenvalues
can be $\ge 1$, as otherwise the corresponding column of ${\Xi}(\tau)$ could be divided by $J(\tau)-1728$ or
$J(\tau)$, retaining holomorphicity in $\bbH$ but spanning over $\bbC[J]$ a strictly larger space
of weakly holomorphic vvmf. Since by \eqref{A2A3} $\cA_2+\cA_3=\Lambda$,  the trace of $\Lambda$
is the sum of the eigenvalues of $\cA_2$ and $\cA_3$, and we thus  obtain \eqref{cvalue}.
\qed%\quad\textbf{QED} to Theorem 3.3
\medskip

Theorem 3.3 first appeared in \cite{bib:BG2}, though for restricted $(\rho,w)$. It is generalised
to arbitrary $T$ and arbitrary genus-0 groups in \cite{bib:BG3}.

Theorem 3.3(c) and \eqref{levelt} say that 
\begin{equation}\label{xiexp2}\Xi(\tau)=\mathrm{j}_2(w;\tau)^{-1}P_2^{-1}q_2^{\lambda_2}
\sum_{n=0}^\infty{\Xi}_{[n]}q_2^n=\mathrm{j}_3(w;\tau)^{-1}P_3^{-1}q_3^{\lambda_3}
\sum_{n=0}^\infty{\Xi}_{\langle n\rangle}q_3^n\end{equation}
where ${\Xi}_{[0]}$ and ${\Xi}_{\langle 0\rangle}$ are invertible, $\lambda_2,
\lambda_3$ are diagonal, $\lambda_2$ has $\alpha_i$ diagonal values equal to $i/2$
for $i=0,1$, $\lambda_3$ has $\beta_j$ diagonal entries equal to $j/3$ for $j=0,1,2$,
$P_2SP_2^{-1}=e^{2\pi\Ii\lambda_2}$, and $P_3UP_3^{-1}=e^{2\pi\Ii\lambda_3}$.
The key properties here are the bounds $0\le (\lambda_2)_{ii}<1$ and $0\le (\lambda_3)_{jj}<1$,
and the invertibility of   ${\Xi}_{[0]}$ and ${\Xi}_{\langle 0\rangle}$.

Given $\Lambda$ and $\chi$, it is easy to solve \eqref{fundeq1} recursively: 
\begin{equation}\label{eq:recurs}
[\Lambda_w,{\Xi}_{(n)}]+n{\Xi}_{(n)}=\sum_{l=0}^{n-1}{\Xi}_{(l)}
\left(f_{n-l}\FB_w+\frac{w}{12}t_{n-l}+g_{n-l}(\chi_w+[\FB_w,\chi_w])\right)
\end{equation}
for $n\ge 2$, where
we write $E_2(\tau)=\sum_{n=0}^\infty t_nq^n=1-24q-\cdots$, $(J(\tau)-984)\Delta(\tau)/E_{10}(\tau)=\sum_{n=0}^\infty f_{n}q^n=1+0q+\cdots$ and
$\Delta(\tau)/E_{10}(\tau)=\sum_{n=0}^\infty g_{n}q^n=q+\cdots$. 
We require ${\Xi}_{(0)}=1_d$. Note that the $ij$-entry on the left-side of \eqref{eq:recurs} is
$\left(\FB_{ii}-\FB_{jj}+n\right){\Xi}_{(n)\, ij}$, so  \eqref{eq:recurs} allows us to recursively
identify all  entries of ${\Xi}_{(n)}$, at least when all 
$|\FB_{jj}-\FB_{ii}|\ne n$.   Indeed, 
 $\FB_{jj}-\FB_{ii}$ can never lie in
$\bbZ_{\ge 2}$, thanks to this recursion, since then the value of ${\Xi}_{(n)\, ij}$ would 
be unconstrained, contradicting uniqueness of the solution to \eqref{fundeq1} with
${\Xi}_{(0)}=1_d$.

Theorem 3.3 tells us that $\bbX(\tau)\in\cM^!_w(\rho)$ iff $\bbX(\tau)={\Xi}(\tau)P(J)$, where $P(J)\in\bbC^d[J]$.
The basis  $\bbX^{(j;n)}(\tau)$ of \eqref{eq:basis} can be easily found recursively from this \cite{bib:BG2}. They can also be found as follows.  Define the 
generating function\begin{equation}
\mathcal{X}_{ij}(\tau,\sigma):=q^{-\FB_{ii}-1}\sum_{n=0}^{\infty}\left[\Fa^{(j;n)}\left(\tau\right)
\right]_{i}z^{n}=\frac{\delta_{ij}}{q-z}+\sum_{m=1}^\infty\sum_{n=0}^\infty\mathcal{X}_{(m)\,i}^{(j;n)}q^{m-1}
z^{n}\,,\label{eq:gfdef}\end{equation}
where we write $z=e^{2\pi\Ii\sigma}$. Then, writing $J'=\Der_0 J=-E^2_4E_6/\Delta$, we have \cite{bib:BG2}
 \begin{equation}
\mathcal{X}\left(\tau,\sigma\right)=\frac{{J'}\left(\sigma\right)\,
q^{-\FB-1_d}}{J\!\left(\tau\right)-J\!\left(\sigma\right)}\Xi(\tau)\,\Xi({\sigma})^{-1}z^{\FB}\:.
\label{eq:wrec4}\end{equation} 

We call ${\Xi}(\tau)$ in Theorem 3.3(b) the \textit{fundamental matrix} associated to $\Lambda$. A 
$d\times d$ matrix $\Xi(\tau)$ is a fundamental matrix for $(\rho,w)$ iff all
columns lie in $\cM^!_w(\rho)$, and $\Xi(\tau)=q^\Lambda(1_d+\sum_{n=1}^\infty\Xi_{(n)}q^n)$
where Tr$\,\Lambda=c_{(\rho,w)}$. The reason is that $\cP_\Lambda:\cM_w^!(\rho)\rightarrow
\bbC^d[q^{-1}]$ is then surjective, so by the index formula it must also be injective.

The determinant of any fundamental matrix is now easy to compute \cite{bib:BG2}:
\begin{equation}\label{eq:detxi}\mathrm{det}\,\Xi(\tau)=E_4(\tau)^{\beta_1+2\beta_2}
 E_6(\tau)^{\alpha_1}\Delta(\tau)^{(dw-4\beta_1-8\beta_2-6\alpha_1)/12}\,,\end{equation}
 where $\alpha_i=\alpha_i(\rho,w),\beta_j=\beta_j(\rho,w)$ are the eigenvalue multiplicities of Section 3.1. Indeed, the
 determinant is a scalar modular form (with multiplier); use $E_4(\tau)$
 and $E_6(\tau)$ to factor off the zeros at the elliptic points (which we can read off from \eqref{xiexp2}), and note that the resulting modular
 form has no zeros in $\bbH$ and hence must be a power of $\Delta(\tau)$, where the
 power is determined by the weight.
 
A very practical way to obtain bijective $\Lambda$ and $\chi$, and hence a fundamental 
matrix ${\Xi}(\tau)$, is through cyclicity:

\medskip\noindent\textbf{Proposition 3.2.}  \textit{Suppose $(\rho,w)$ is admissible,
and $T$ is diagonal. Suppose $\bbX(\tau)\in\cM^!_w(\rho)$, and the components of $\bbX(\tau)$ are
linearly independent over $\bbC$. Then $\cM^!_w(\rho)=\bbC[J,\nabla_{1,w},\nabla_{2,w},\nabla_{3,w}]\bbX(\tau)$.}\medskip

\noindent\textit{Proof.} Let $\cM_\bbX:=\bbC[J,\nabla_{1,w}]\bbX(\tau)$. Since $\cM_\bbX$ is a module
of a PID $\bbC[J]$, it is a sum of cyclic submodules $\bbC[J]\bbY^{(i)}(\tau)$. Each $\bbC[J]\bbY^{(i)}(\tau)$ is 
torsion-free (by looking at leading powers of $q$). So $\cM_\bbX$ must be free of some rank $d'$.
Because it is a submodule of the rank $d$ module $\cM_w^!(\rho)$ (and again using the
fact that $\bbC[J]$ is a PID), $d'\le d$. That $d'=d$ follows by computing the determinant of
the $d\times d$ matrix with columns $\nabla_{1,w}^{i-1}\bbX(\tau)$: that determinant equals $(2\pi\Ii)^{1-d}
(E_{4}E_6/\Delta)^{(d-1)d/2}$ times the Wronskian of $\bbX(\tau)$, which is nonzero by Lemma 2.2(a). 

Let $\Xi_\bbX(\tau)$ be the matrix formed by those $d$ generators of $\cM_\bbX$. Because 
$\nabla_{1,w}\cM_\bbX\subseteq\cM_\bbX$, the argument of Theorem 3.3 applies and
$\Xi_\bbX(\tau)$ satisfies analogues of \eqref{fundeq1} and hence \eqref{fundeq2}. This means $\cM_\bbX$
will have its own analogues $\Lambda_\bbX,\cA_{\bbX\, 2},\cA_{\bbX\,3}$ (their exponentials
$e^{2\pi\Ii\Lambda_\bbX}$ etc will be conjugate to $e^{2\pi\Ii\Lambda}$ etc). The trace of
$\Lambda_\bbX$ will equal the trace of $\cA_{\bbX\, 2}+\cA_{\bbX\,3}$,
for the same reason it did in Theorem 3.3. Now, $\cP_{\Lambda_\bbX}$
is an isomorphism when restricted to $\cM_\bbX$, so when extended to $\cM_w^!(\rho)$ will also
have trivial cokernel. Thus the dimension of $\cM_w^!(\rho)/\cM_\bbX$ will equal the dimension
of ker$\,\cP_{\Lambda_\bbX}$ which, by the index formula \eqref{index}, equals 
$$\mathrm{Tr}\,\Lambda_\bbX-\mathrm{Tr}\,\Lambda=(\mathrm{Tr}\,\cA_{\bbX\,2}-\mathrm{Tr}\,\cA_2)+(\mathrm{Tr}\,
\cA_{\bbX\,3}-\cA_3)\,.$$

This means that if $\cM_\bbX\ne \cM^!_w(\rho)$, then at least one eigenvalue of $\cA_{\bbX\,2}$
or $\cA_{\bbX\,3}$ is $\ge 1$.  Suppose one of $\cA_{\bbX\,2}$ is.  Then by \eqref{xiexp2} every
component of some row of
$\Xi_\bbX(\tau)$ has order $\ge 2$ at $\tau=\Ii$. This means every vvmf in $\cM_\bbX$ has some
component with a zero at  $\tau=\Ii$ of order $\ge 2$. Hit each column of $\Xi_\bbX(\tau)$ with
$\nabla_{2,w}=\Delta^{-1}E_4^2\Der_w^2$: it will reduce the order of that zero everywhere by 2
(because at $\tau=\Ii$ $\nabla_{2,w}$ looks like $a\frac{\mathrm{d}^2}{\mathrm{d}\tau_2^2}
+b\frac{\mathrm{d}}{\mathrm{d}\tau_2} +c$ for $a\ne 0$).
This means the module generated over $\bbC[J,\nabla_{1,w}]$ by the columns of
$\Xi_\bbX(\tau)$ and $\nabla_{2,w}\Xi_\bbX(\tau)$ is strictly greater than $\cM_\bbX$, as some
vvmf in it has a smaller order at that component than all vvmf in $\cM_\bbX$. So repeat
this argument with $\cM_\bbX$ replaced with this extension. If instead $\cA_{\bbX\,3}$
has an eigenvalue $\ge 1$, then some component of every vvmf in $\cM_\bbX$ has a zero 
at $\tau=\xi_6$  of order
$\ge 3$, so use $\nabla_{3,w}$, which reduces the order of that zero by 3. Eventually
all eigenvalues will be $<1$, in which case the dimension of $\cM^!_w(\rho)/\cM$ will be 0. \qed%\qquad\textbf{QED} to Proposition 3.2
\medskip

To indicate the nontriviality of our theory, we get a 1-line proof of the solution \eqref{eq:inequal}
to the Deligne-Simpson problem for $\overline{\Gamma}$, at least for most $\rho$ (the general
case requires slightly more work). Let $\rho$ be any $\overline{\Gamma}$-representation with
$T$ diagonal, and let $\Lambda$ be bijective for $(\rho,0)$ and $\Xi(\tau)$ a corresponding
fundamental matrix. Then as long as all $\Lambda_{ii}\ne 0$, the columns of the derivative
$\Der_0\Xi$ will span a free rank-$d$ submodule of $\cM^!_2(\rho)$ over $\bbC[J]$,
on which $\cP_\Lambda$ is surjective. Thus $c_{(\rho,0)}\le c_{(\rho,2)}$ and so by
\eqref{tracek} we obtain $\alpha_1\ge\beta_0$. (It is clear that things are more subtle when
some $\Lambda_{ii}=0$, as this inequality fails for $\rho=1$!) The other 
inequalities $\alpha_i\ge \beta_j$ follow by comparing $c_{(\rho,2k)}$ and $c_{(\rho,2k+2)}$
in the identical way.

As mentioned earlier, $S$ and $U$ are conjugate to $e^{2\pi\Ii\cA_2}$
and  $e^{2\pi\Ii\cA_3}$, respectively, but identifying precisely \textit{which} conjugate is
a transcendental and subtle question. For example, we 
see in Section 4.2 below that when $d=2$, they are related by Gamma function values.

\subsection{Holomorphic vvmf}

Until this point in the paper, our focus has been on weakly holomorphic vvmf of fixed weight,
i.e. vvmf holomorphic everywhere except at the $\Gamma$-orbit of $\Ii\infty$.
The reason is that structurally it is the simplest and most fundamental. For example, it is acted on by the
ring of (scalar) modular functions holomorphic away from $\Gamma\Ii\infty$, which for any genus-0 Fuchsian group is  a PID.
By contrast, the holomorphic vvmf (say of arbitrary even integral weight) is a module over the ring of holomorphic modular
forms, which in genus-0 is usually not even polynomial.
However, there is probably more interest in holomorphic vvmf, so it is to these we now turn. 
The two main questions we address are the algebraic structure (see Theorem 3.4
below), and dimensions (see Theorem 3.5 next subsection).

\medskip\noindent\textbf{Definition 3.1.} \textit{Let $(\rho,w)$ be admissible, 
$T$  diagonal, and $\lambda$ any exponent. Define 
\begin{equation}\cM^\lambda_w(\rho):=\mathrm{ker}\,\cP_{\lambda-1_d}=\left\{\bbX(\tau)\in\cM^!_w(\rho)\,\,\,|\,\,\,
\bbX(\tau)=q^\lambda\sum_{n=0}^\infty \bbX_{(n)}q^n \right\}\end{equation}   
and $\cM^\lambda(\rho)=\coprodd_{k\in\bbZ}\cM^\lambda_{w+2k}(\rho)$. We call any $\bbX(\tau
)\in
\cM^\lambda(\rho)$,} $\lambda$-holomorphic.\medskip

For example, for the trivial representation, $\cM^{0}(1)$ are the modular forms $\frakm=\bbC[E_4,E_6]$,
while $\cM^{1}(1)$ are the cusp forms $\frakm\Delta(\tau)$. More generally, define
 $\lambda^{hol}$ to be the unique exponent with $0\le \mathrm{Re}\,\lambda_{ii}<1$
 for all $i$. Then $\cM^{\lambda^{hol}}(\rho)$ coincides with the usual definition of holomorphic
 vvmf. Choosing $0<\mathrm{Re}\,\lambda_{ii}\le 1$ would give the vector-valued cusp forms.
 
 Theorem 3.2 computes dim$\,\cM_{w+2k}^\lambda(\rho)$ for all sufficiently
 large $|k|$:
 
 \medskip\noindent\textbf{Lemma 3.2.} \textit{Let $(\rho,w)$ be admissible and $T$ diagonal.
 Let $\lambda$ be any exponent, and $\Lambda$ any bijective exponent.}
 
 \smallskip\noindent\textbf{(a)} \textit{For any $k\in\bbZ$, $\cM^\lambda_{w+2k}(\rho)$ is finite-dimensional,
and obeys the bound}
\begin{equation}\label{holbound}
\mathrm{dim}\,\cM^\lambda_{w+2k}(\rho)\ge\mathrm{max}\left\{0,\frac{w+2k+2}{12}d+\frac{\alpha_k}{2}
+\frac{\beta_k-\beta_{k+2}}{3}-\mathrm{Tr}\,\lambda\right\}\,.\end{equation}
 
  \smallskip\noindent\textbf{(b)} \textit{Choose any $m,n\in\bbZ$ satisfying
  $\mathrm{ Re}\,(\Lambda+m1_d)\ge \mathrm{Re}\,\lambda\ge\mathrm{Re}\,(\Lambda-n1_d)$ entrywise.
Then   $\cM^\lambda_{w-2k}(\rho)=0$ when $k=6n+6$ or $k\ge 6n+8$, and equality holds
in \eqref{holbound}  whenever $k=6m-6$ or $k\ge 6m-4$.}   

\smallskip\noindent\textbf{(c)} \textit{Let  $w_0$ be the weight with smallest real part for which
$\cM^{\lambda}_{w_0}(\rho)\ne 0$ and write $\epsilon=0,1$ for $d$ even, odd respectively. Suppose  $\rho$ is irreducible and not 1-dimensional. Then $w_0$ 
 satisfies the bounds}
\begin{equation}\label{eq:bounds}
\frac{12}{d} \mathrm{Tr}\,\lambda+1-d\le w_0\le \frac{12}{d}\mathrm{Tr}\,\lambda-\frac{3\epsilon}{d}
\,.\end{equation}

\medskip\noindent\textit{Proof.}  Theorem 3.2(a) gives finite-dimensionality. The bound \eqref{holbound}
follows from $\mathrm{dim\,ker}\,\cP_{\lambda-1_d}\ge \mathrm{index}\,\cP_{\lambda-1_d}$,
the index formula in Theorem 3.2(a), and \eqref{alphabet}. Note that $\Delta(\tau)^k\cM^!_w(\rho)=\cM^!_{w+2k}(\rho)$
for any $k\in\bbZ$ so $\Lambda+k1_d$ is bijective for $(\rho,w+12k)$. Hence
$\cP_{\lambda-1_d}$ is injective on $\cM^!_{w-12(n+1)}(\rho)$ because $\cP_{\Lambda
-(n+1)1_d}$ is, while $\cP_{\lambda-1_d}$ is surjective on $\cM^!_{w+12(m-1)}(\rho)$
because $\cP_{\Lambda+(m-1)1_d}$ is. This proves (b) for those weights.
Now, for any $k\ge 2$ there is a scalar  modular form
$f(\tau)\in\frakm$ of weight $2k$ with nonzero constant term, so $f(\tau)^{-1}\in\bbC[[q]]$ and the 
surjectivity of $\cP_{\lambda-1_d}$ on $\cM_{w+12(m-1)}^!(\rho)$ implies the surjectivity on
$\cM^!_{w+12(m-1)+2k}(\rho)$. More directly, injectivity of $\cP_{\lambda-1_d}$ on $\cM_{w-12(n+1)}^!(\rho)$ implies the injectivity on
$\cM^!_{w-12(n+1)+2k}(\rho)$.

Now turn to (c). 
We know that dim$\,\mathcal{M}_{w}^{\lambda}(\rho)>0$ for any $w$ for which Tr$\,\FB>\mathrm{
Tr}\,\lambda-d$, since dim$\,\mathrm{ker}\,\cP_{\lambda-1_d}\ge\mathrm{Index}\,\cP_{\lambda-1_d}>0$.
 Using \eqref{eq:Trbound}, we obtain the upper bound of  \eqref{eq:bounds}.

Choose any nonzero $\FA\in\mathcal{M}^{\lambda}_{w}(\rho)$. Then its components must
be linearly independent over $\bbC$, because they span a subrepresentation of the
irreducible $\rho$. Therefore, Lemma 2.2(a) says Wr$(\bbX)(\tau)\in\cM^!_{d(w+d-1)}(\mathrm{det}
\,\rho)$ is nonzero, with leading power
of $q$ in Tr$\,\lambda+\bbZ_{\ge 0}$. Hence  Wr$(\bbX)\tau/\Delta^{\mathrm{Tr}\,\lambda}
(\tau)$ lies in
$\mathcal{M}_{d(w+d-1)-12\mathrm{Tr}\,\lambda}^{0}(\varmu_u)$ for some $u\in2\bbZ$, from which follows
 the lower bound of \eqref{eq:bounds}. \qed%\quad \textbf{QED} to Lemma 3.2
\medskip 

The lower bound in \eqref{eq:bounds} is due to Mason \cite{bib:Mas} (he proved it for
$\lambda=\lambda^{hol}$ but the generalisation given here is trivial). The $d\ne 1$ assumption
in (c) is only needed for the upper bound.

Theorem 3.3(a) tells us the space $\cM^!_w(\rho)$ of weakly holomorphic vvmf is a free
module  of rank $d$ over $\bbC[J]$. The analogous statement for holomorphic vvmf,
namely that $\cM^{\lambda^{hol}}(\rho)$ is free of rank $d$ over $\frakm$, is implicit in
\cite{bib:EZ} (see the Remark there on page 98). It was also proved independently in
\cite{bib:MM}, and independently but simultaneously we obtained the following generalisation.
The proof of freeness given here  is far simpler than in \cite{bib:MM}, is more general (as it applies to arbitrary $\lambda$),
gives more information (see (b) below), and generalises directly to arbitrary $T$ and arbitrary 
genus-0 groups \cite{bib:Ga}. 

\medskip\noindent\textbf{Theorem 3.4.} \textit{Let $(\rho,w)$ be admissible and $T$
diagonal. Choose any exponent $\lambda$. Let $\alpha_j=\alpha_j(\rho,w)$ and $\beta_j=
\beta_j(\rho,w)$.}

\smallskip\noindent\textbf{(a)} \textit{$\cM^\lambda(\rho)$ is a free module
over $\frakm=\bbC[E_4,E_6]$, of rank $d$.}

\smallskip\noindent\textbf{(b)} 
\textit{Let $w_0=w^{(1)}\le w^{(2)}\le\cdots\le w^{(d)}$ be the weights of the free generators.
Then precisely $\alpha_i$ of the $w^{(j)}$ are congruent mod 4 to $w+2i$, and
precisely $\beta_i$ of them are congruent mod 6 to $w-2i$. Moreover, $\sum_j w^{(j)}
= 12\,\mathrm{Tr}\,\lambda$. Let ${\Xi}^{\lambda}(\tau)$ be the matrix 
obtained by putting these $d$ generators  into $d$ 
columns. Then (up to an irrelevant nonzero constant)}
 \begin{equation}\nonumber\mathrm{det}\,{\Xi}^{\lambda}(\tau)=\Delta^{\mathrm{Tr}\,\lambda}(\tau)\,.\end{equation}

%\smallskip\noindent\textbf{(c)} \textit{As a module over $\bbC[E_4,E_6,\Der]$,
%precisely max$\{n_0,n_1\}$ generators are needed, namely $n_0$ 
%in $\mathcal{M}_0$ and max$\{0,\alpha_1-\beta_0\}$ in $\mathcal{M}_1$. (restriction on
%whether 1 is summand?)}

\medskip\noindent\textbf{Proof.} 
Let $w_0\in w+2\bbZ$ be the weight of minimal real part with $\cM^\lambda_{w_0}(\rho)\ne 0$
--- this exists by Lemma 3.2(b).  Write $\mathcal{M}=
\mathcal{M}^{\lambda}(\rho)$, $w_k=w_0+2k$ and $\mathcal{M}_l
=\mathcal{M}^{\lambda}_{w_l}(\rho)$. For $\bbX(\tau)\in\cM_w(\rho)$, recall from 
\eqref{expansions} that the constant term $\bbX_{[0]}$ at $\tau=\Ii$ is $\bbX(\Ii)/E_4(\Ii)^{w/4}$. 
Fix $S_0=e^{\pi\Ii w_0/2}S$;
then $e^{\pi\Ii w_k/2}S=(-1)^kS_0$ and so for any $\bbX\in\cM_k$, its constant term satisfies
$S_0\bbX_{[0]}=(-1)^k\bbX_{[0]}$ thanks to \eqref{SUtrans}.

Find $\Fa^{(i)}(\tau)\in\mathcal{M}_{l_i}$ with the property that, for any $k\ge 0$,
the space of constant terms $\Fa_{[0]}$, as $\Fa(\tau)$ runs over all $\cup_{l=0}^k
\mathcal{M}_l$, has a basis given by the constant terms $\Fa^{(i)}_{[0]}$ for
those $\Fa^{(i)}(\tau)$ in $\cup_{l=0}^k\mathcal{M}_l$ (i.e. for those $i$ with $l_i\le k$). This is done
recursively with $k$. We will show that these $\Fa^{(i)}(\tau)$ freely generate $\cM$.

The key observation is the following. 
Consider any $\FA(\tau)\in\mathcal{M}_k$. Then by definition of the $\bbX^{(i)}(\tau)$,
$\Fa_{[0]}=\sum_ic_i \Fa^{(i)}_{[0]}$ where $c_i=0$ unless $l_i\le k$ and  
$l_i\equiv k$ (mod 2). Now observe that  the constant term $\bbX'_{[0]}$ of
$\Fa'(\tau)=\FA(\tau)-\sum_ic_i \,E_4(\tau)^{(k-l_i)/2}
\Fa^{(i)}(\tau)$ is 0, so $\Fa'(\tau)/E_6(\tau)\in\mathcal{M}_{k-3}$.

One consequence of this observation is that, by an easy induction on $k$,
any $\Fa(\tau)\in\cM$ must lie in $\sum_i\frakm\Fa^{(i)}(\tau)$. Another consequence is
that there are exactly $d$ of these $\Fa^{(i)}(\tau)$,  in particular their constant terms
form a basis for $\bbC^d$. To see this, take any fundamental matrix ${\Xi}(\tau)$ at weight $w_0$.
Then for sufficiently large $l$, each column of $\Delta(\tau)^l{\Xi}(\tau)$ is in 
$\mathcal{M}_{12l}$.
The constant terms of  $\Delta(\tau)^l{\Xi}(\tau)$ (which have $S_0$-eigenvalues +1) and of 
$\Delta(\tau)^l\Der_{w_0}{\Xi}(\tau)$ (which have $S_0$-eigenvalues $-1$) must have rank
$\alpha_0(\rho,w_0)$ and $\alpha_1(\rho,w_0)$, respectively, as otherwise the columns of $\Delta(\tau)^l{\Xi}(\tau)$ 
would be linearly dependent over $\frakm$, contradicting their linear independence
over $\bbC[J]$. This means of course that exactly $\alpha_0(\rho,w_0)$ of these $\bbX^{(i)}(\tau)$
 have $l_i$ even, and exactly $\alpha_1(\rho,w_0)$ have $l_i$ odd.

Thus these $d$ $\Fa^{(i)}(\tau)$ generate over $\frakm$ all of $\mathcal{M}$.
To see they are  linearly independent over $\frakm$, suppose we have a relation
 $\sum_i p_i(\tau)\,\Fa^{(i)}(\tau)=0$, 
for modular forms  $p_i(\tau)\in\frakm$ which do not share a (nontrivial) common divisor. The constant
term at $\Ii$ of that relation  reads $\sum_i p_i(\Ii)\,\Fa^{(i)}_{[0]}=0$  and hence each 
$p_i(\Ii)=0$, since the $\bbX^{(i)}_{[0]}$ are linearly independent by construction. This forces all  $p_i(\tau)$ to be 0, since otherwise we could divide them all by  $E_6$, which would contradict the hypothesis
that they share no common divisor.

The identical argument applies to the constant terms $\bbX_{\langle 0\rangle}=\bbX(\xi_6)/E_6
(\xi_6)^{w/6}$ at $\tau=\xi_6$; this implies
that exactly $\beta_i(\rho,w_0)$ generators $\bbX^{(j)}(\tau)$ have $l_j\equiv i$ (mod 3).

Form the $d\times d$ matrix $\Xi^\lambda(\tau)$ from these $d$ generators $\bbX^{(i)}(\tau)$ and call the determinant $\delta(\tau)$.
The linear independence of the
constant terms of the generators at the elliptic fixed points, says that $\delta(\tau)$
 cannot vanish at any
elliptic fixed point. $\delta(\tau)$ also can't have a zero anywhere else in $\uhp$. To see this, first
note that a zero at $\tau^*\in\bbH$ implies there is nonzero row vector ${v}\in\bbC^d$ such that
$v\Xi^\lambda(\tau)=0$ and hence $v\bbX(\tau^*)=0$ for any $\bbX(\tau)\in\mathcal{M}^{\lambda}(\rho)$. But  \eqref{eq:detxi} says the
determinant of any fundamental matrix $\Xi(\tau)$ for $(\rho,w_0)$ can only vanish at elliptic fixed points and cusps, and so
 $v{\Xi}(\tau^*)\ne 0$ for any fundamental matrix, and hence $v\bbY
(\tau^*)\ne 0$ for some $\bbY(\tau)\in\cM^!_{w_0}(\rho)$. To get a contradiction, choose $N$
big enough so that $\Delta(\tau)^N\bbY(\tau)\in\cM$. This means $\delta(\tau)$ is a scalar modular form (with 
multiplier) which doesn't vanish anywhere in $\bbH$. Hence $\delta(\tau)$ must be a power of $\Delta(\tau)$,
and considering weights we see this must be $\delta(\tau)=\Delta^{\sum_iw_{l_i}/12}(\tau)$. We compute that
sum over $i$, shortly.

Find the smallest $\ell$ such that $c_{(\rho,w_0+2\ell)}>\mathrm{Tr}
\,\lambda$ (recall \eqref{tracek}) and put $w_0':=w_\ell$.  
Define $n_k'=0$ for $k<0$, and \begin{equation}
n'_k=\frac{w_0'+2k+2}{12}d+\frac{\alpha_k}{2}+
\frac{\beta_k-\beta_{2+k}}{3}-\textrm{Tr}\,\lambda
\end{equation}
 for $k\ge 0$. Using \eqref{alphabet}, the numbers $n'_k$ are the values max$\{c_{(\rho,w_0'+2k)},0\}$.
From Lemma 3.2 we obtain the equality $n_k=n'_{k-\ell}$
 for $|k|$ sufficiently large.
The \textit{tight Hilbert--Poincar\'e series} $H^{\lambda}_{tt}(\cM;x):=\sum_k n'_kx^{w_{k+\ell}}$
equals
\begin{equation}\label{eq:fake}
H^\lambda_{tt}(\cM;x)=x^{w_0'}\frac{n'_0+n'_1x^2+(n'_2-n'_0)x^4+(n'_3-n'_1-n'_0)x^6
+(n'_4-n'_2-n'_1)x^8}{(1-x^4)(1-x^6)}\,,
\end{equation}
by a simple calculation (the significance of $H_{tt}^\lambda(x)$ is explained in Proposition
3.3 below). From the numerator we read off the weights of the hypothetical `tight' generators: write $w'{}^{(i)}=w_0'$ for $1\le i\le n_0'$, $ \ldots$, $w'{}^{(i)}=w_0'+8$ for $d-n_4'+n_2'+n_1'<i\le d$.
We know that the actual Hilbert--Poincar\'e series, $H^\lambda(\cM;x)=\sum_kn_kx^{w_k}$, minus the tight one, equals a finite sum of
terms, since $n_k=n'_{k-\ell}$ for $|k|$ large. Therefore $(1-x^4)(1-x^6)(H^{\lambda}(x)-H^{\lambda}_{tt}(x))=\sum x^{w^{(i)}}-\sum x^{w'{}^{(i)}}$
is simply a polynomial identity.
Differentiating with respect to $x$ and setting $x=1$ gives $\sum w^{(i)}=\sum w'{}^{(i)}$,
and the latter is readily computed to be $12\,\mathrm{Tr}\,\lambda$. 
\qed%\quad \textbf{QED} to Theorem 3.4
\medskip

This freeness doesn't seem directly related to that of Theorem 3.3(a).
The reason it is natural to look at local expansions about $\tau=\Ii$ and $\tau=\zeta_6$ in the
proof of part (a) 
is because if $\bbX(\tau)$ is holomorphic and $\bbX(\Ii)=0$, then $\bbX(\tau)/E_6(\tau)$ is also holomorphic
(similarly for $\bbX(\zeta_6)=0$ and $\bbX(\tau)/E_4(\tau)$).

Call an admissible $(\rho,w)$ \textit{tight} if it has the property that for all $k\in\bbZ$,
an exponent $\lambda$ is bijective for $w+2k$ iff Tr$\,\lambda=c_{(\rho,w+2k)}$. In this case we also
say $\rho$ is tight. Generic $\rho$ are tight. We learn in Theorem 4.1 that all irreducible $\rho$
in dimension $d<6$ are tight. For tight $\rho$, most quantities can be easily determined:

\medskip\noindent\textbf{Proposition 3.3.} \textit{Suppose $(\rho,w)$ is admissible and tight,
and $T$ is diagonal. }

\smallskip\noindent\textbf{(a)} \textit{Then $\varmu_u\otimes \rho$ is also tight for any $u\in\bbC$,
as is the contragredient $\rho^*=(\rho^t)^{-1}$.}\smallskip

\noindent\textbf{(b)} \textit{For any exponent $\lambda$ (not necessarily bijective),
either ker$\,\cP_\lambda=0$ (if Tr$\,\lambda\ge c_{(\rho,w)}$) or coker$\,\cP_\lambda=0$
(if Tr$\,\lambda\le c_{(\rho,w)}$). Moreover,
\begin{equation}\mathrm{dim}\,\cM^\lambda_w(\rho)=\mathrm{max}\{0,c_{(\rho,w)}+d-\mathrm{Tr}
\,\lambda\}\,,\end{equation}
and the Hilbert--Poincar\'e series $H^\lambda(x)$ of $\cM^\lambda(\rho)$ equals the tight Hilbert--Poincar\'e
series $H_{tt}^\lambda(x)$ of \eqref{eq:fake}.}\medskip

The proof of (a) uses the equality $\cM^!_{w+u}(\varmu_u\otimes \rho)=\Delta^{u/12}(\tau)
\cM^!_w(\rho)$, as well as the duality in Proposition 3.4 below. To prove (b), let $\Lambda$ be
the unique bijective exponent matrix for $(\rho,w)$ 
satisfying $\Lambda_{ii}=\lambda_{ii}-1$ for all $i\ne 1$. Write $n=\Lambda_{11}+1-\lambda_{11}=c_{(\rho,w)}+d-
\textrm{Tr}\,\lambda\in\bbZ$. If $n\ge 0$ $\cP_\lambda$ inherits the surjectivity of $\cP_\Lambda$,
while if $n\le 0$ it inherits the injectivity. The index formula \eqref{index} gives the dimension,
 which by \eqref{alphabet} equals $n_{k-\ell}'$. This is why $H_{tt}^\lambda$ arises.

As long as $(\rho,w)$ is tight, this argument tells us how to  find a basis for $\mathcal{M}^{\lambda}_{w}(\rho)$. Let $\Lambda$ and $n$ be as above.  For $n>0$ a basis for 
$\mathcal{M}^{\lambda}_{w}(\rho)$ consists of the basis vectors
$\bbX^{(1;i)}$  (see Section 3.2) for $0\le i< n$. Incidentally, the name `tight' refers to the fact that 
the numerator of $H_{tt}^\lambda$ is maximally bundled together.

There are other constraints on the possible weights $w^{(i)}$ of Theorem 3.4(b). A useful
observation in practice is that if $\rho$ is irreducible, then the set $\{w^{(1)},\ldots,w^{(d)}\}$
can't have gaps, i.e. for $n=(w^{(d)}-w^{(1)})/2$,
\begin{equation}\label{eq:gap}
\{w^{(1)},\ldots,w^{(d)}\}=\{w^{(1)},w^{(1)}+2,w^{(1)}+4,\ldots,w^{(1)}+2n\}\,.\end{equation}
The reason is that when $w^{(1)}+2k$ doesn't equal any $w^{(i)}$, then when $w^{(j)}<w^{(1)}+2k$,
$\Der^{k-(w^{(j)}-w^{(1)})/2}\bbX^{(j)}(\tau)=\sum_{w^{(l)}<w^{(1)}+2k}f_{jl}(\tau)\bbX^{(l)}(\tau)$ for $f_{jl}(\tau)\in\frakm$, where the sum
is over all $l$ with $w^{(l)}<w^{(1)}+2k$. If in addition $w^{(1)}+2k<w^{(d)}$ (i.e. $w^{(1)}+2k$
is a \textit{gap}), then the Wronskian Wr$(\bbX^{(1)})$ would have to vanish, contradicting
irreducibility.

Several papers (e.g. \cite{bib:Mas,bib:Mar}) consider a `cyclic' class of vvmf where the components
of $\bbX(\tau)$ span the solution space to a monic modular differential equation $\Der_k^d+
\sum_{l=0}^{d-1} f_l(\tau)D_k^l=0$, where each $f_l(\tau)\in\frakm$ is of weight $2d-2l$. In this
accessible case the free generators $\bbX^{(i)}(\tau)$ can be taken to be $\bbX(\tau),\Der_k\bbX(\tau),
\ldots,\Der_k^{d-1}\bbX(\tau)$, i.e. the weights are $\alpha^{(i)}=\alpha^{(1)}+(i-1)2$. This means 
the corresponding $\rho$ \textit{cannot} be tight, when $d\ge 6$. Indeed, the multipliers of such vvmf are exceptional, requiring the
multiplicities of $\rho$ to satisfy $|\alpha_i-\alpha_j|\le 1$ and $|\beta_i-\beta_j|\le 1$.
Recall that the connected components of the moduli space of $\overline{\bG}$-representations
are parametrised by these multiplicities; these particular components are of maximal dimension.
For example, when $d=6$, such a representation $\rho$ can lie in only 1 of the 12 possible
connected components, and these $\rho$ define a 6-dimensional subspace inside that
7-dimensional component.

\subsection{Serre duality and the dimension formula}

A crucial symmetry of the theory is called the \textit{adjoint} in the language of Fuchsian equations, 
and shortly we reinterpret this as Serre duality.

\medskip\noindent
\textbf{Proposition 3.4.} \textit{Let $(\rho,w)$ be admissible, $T$ diagonal, and $\Lambda$ 
bijective, and let ${\Xi}(\tau)$ be the associated  fundamental matrix of $\cM_{w}^!(\rho)$. 
Let $\rho^*$ denote the
{contragredient} $({\rho}^{-1})^t$ of $\rho$.
Then $\cM_{w}^!(\rho)$
and $\cM_{2-w}^!(\rho^*)$ are naturally isomorphic as $\bbC[J,\nabla_1,\nabla_2,\nabla_3]$-modules.
Moreover,}
\begin{align}%\FB^*_\bac=\; & 1_d-\FB_\bac\,,\label{eq:starlambda}\\
%{\mathcal{X}}^*=\; & -\,\trans{\mathcal{X}}\:,\label{eq:starX}\\
{\Xi}^*(\tau)=\; & E_4(\tau)^2E_6(\tau)\,\Delta(\tau)^{-1}\left(\Xi({\tau})^{t}\right)^{-1}\,,\label{eq:serre}\\
% Moreover, the matrix-valued generating functions $\frak{X}(q,z)$ obey the duality
{\mathcal{X}}^*(\sigma,\tau)=\;&-{\mathcal{X}}(\tau,\sigma)^t\,,\label{eq:zqsym}\end{align}
\noindent\textit{where $\Xi^*$ is the fundamental matrix of $\cM^!_{2-w}(\rho^*)$ corresponding to
the bijective exponent $\Lambda^*=-1_d-\Lambda$, and $\mathcal{X},\mathcal{X}^*$ are the generating functions
\eqref{eq:gfdef} for $(\rho,w)$ and $(\rho^*,2-w)$ respectively. In other words the $q^{m+\FB_{ii}}$ coefficient $\mathcal{X}^{(j;n-1)}_{(m)\,i}$ of the  basis
vector $\Fa^{(j;n-1)}_{i}(\tau)$ is the negative of the $q^{n+\FB_{jj}^*}$ coefficient $\mathcal{X}^{*\,(i;m-1)}_{(n)\,j}$
of the  basis vector $\Fa_{j}^{*\left(i;m-1\right)}(\tau)$, for all $m,n\ge 1$.}\medskip
%For $m=n=1$ this reduces to \eqref{eq:starX}.

 \noindent\textbf{Proof.}  Define $\Xi^*(\tau)$ by \eqref{eq:serre}; to show it is the fundamental
 matrix associated to the bijective exponent $\Lambda^*=-1_d-\Lambda$, we need to show
that $\Xi^*(\tau)=q^{-1_d-\Lambda} (1_d+\sum_{n=1}^\infty \Xi^*_{(n)}q^n)$ (this is clear),
that $-1_d-\Lambda$ has trace $c_{(\rho^*,2-w)}$ (we'll do this next), and that the columns
of $\Xi^*(\tau)$ are in $\cM^!_{2-w}(\rho^*)$ (we'll do that next paragraph). Using \eqref{alphabet}, we find $\alpha_i(\rho^*,
 2-w)=\alpha_{i+1}(\rho,w)$, and $\beta_j(\rho^*,2-w)=\beta_{2-j}(\rho,w)$, so 
 we compute from \eqref{cvalue} that $c_{(\rho^*,2-w)}=\mathrm{Tr}\,(-1_d-\Lambda)$. 
 
From $\Xi(\gamma\tau)=\widetilde{\rho}_w(\gamma,\tau)\,\Xi(\tau)$, we get $(\Xi(\gamma
 \tau)^t)^{-1}=\widetilde{\rho}^*_{-w}(\gamma,\tau)(\Xi(\tau)^t)^{-1}$. It thus suffices to show 
$\Xi^*(\tau)$ is holomorphic in $\bbH$. We see from \eqref{eq:detxi} that $\Xi(\tau)^{-1}$ is meromorphic everywhere
 in $\bbH^*$, with poles possible only at the elliptic points and the cusp. Locally about $\tau=\Ii$,
 \eqref{xiexp2} tells us 
 \begin{equation}(\Xi(\tau)^t)^{-1}=\mathrm{j}(-w;\tau)\,P_2^tq_2^{-\lambda_2}\left(1_d+\sum_{n=1}^\infty
 ({\Xi}_{[0]}^t)^{-1}{\Xi}_{[n]}^tq^n\right)^{-1}({\Xi}_{[0]}^t)^{-1}\,.\end{equation}
 The series in the middle bracketed factor is invertible at $\tau=\Ii$, because its determinant equals 1 there.
 So every entry of $(\Xi(\tau)^t)^{-1}$ at $\tau=\Ii$ has at worst a simple pole (coming from
 $q_2^{-\lambda_2}=\tau_2^{-2\lambda_2}$). But $J'=-E_4^2E_6/\Delta$ has a simple pole at $\tau=\Ii$.
 Therefore $\Xi^*(\tau)$ is holomorphic at $\tau=\Ii$. Likewise, at $\tau=\xi_6$, the entries
 of $(\Xi(\tau)^t)^{-1}$ have at worse an order 2 pole (coming from $q_3^{-\lambda_3}=\tau_3^{-3
 \lambda_3}$), but $J'$ has an order 2 zero at $\tau=\xi_6$, so $\Xi^*(\tau)$
 is also holomorphic at $\tau=\xi_6$. Thus the columns of $\Xi^*(\tau)$
 lie in $\cM_{2-w}^!(\rho^*)$ (the 2 comes from $J'$). 

 This concludes the proof that $\Xi^*(\tau)$ is a fundamental matrix for $(\rho^*,2-w)$.
  \eqref{eq:zqsym} now follows directly from \eqref{eq:wrec4}.
 \qed%\textbf{QED to Corollary 3.10}
\medskip

A special case of equation \eqref{eq:serre} was found in \cite{bib:BG2}.
An interesting special case of \eqref{eq:zqsym} is that the constant term of any weakly holomorphic
modular form 
$f(\tau)\in\mathcal{M}^!_2(1)$ is 0. To see this, recall first that $(\rho,w)=(1,0)$ has $\Lambda=0$ and
$\Xi(\tau)=\bbX^{(1;0)}(\tau)=1$, i.e. the $q^{m+0}$-coefficient of $\bbX^{(1;0)}_1$ vanishes
for all $m\ge 1$. Then  \eqref{eq:zqsym} says the $q^{1-1}$-coefficient of all $\bbX^{(1;m-1)}_1(\tau)$
must also vanish. Since the $\bbX^{(1;k)}(\tau)$ span (over $\bbC$) all of $\cM^!_2(1)$, the
constant term of any $f(\tau)\in\cM^!_2(1)$ must vanish. The same result holds for any genus-0
group.

The final fundamental ingredient of our theory connects the index formula \eqref{index}
to this duality:

\medskip\noindent\textbf{Theorem 3.5.} \textit{Let $(\rho,w)$ be admissible and $T$ diagonal,
and recall the quantity $c_{(\rho,w)}$ computed in \eqref{cvalue}. Then for any exponent $\lambda$,}
\begin{eqnarray}&\mathrm{coker}\,\cP_{\lambda;(\rho,w)}\cong\left(\cM^{-\lambda}_{2-w}(\rho^*)\right)^*
\,, \label{eq:serre2}\\ & \mathrm{dim}\,\cM^{\lambda}_{w}(\rho)
-\mathrm{dim}\,\cM_{2-w}^{1_d-\lambda}(\rho^*)=c_{(\rho,w)}+d-\mathrm{Tr}\,\lambda\,.\label{eq:index2}\end{eqnarray}

\noindent\textit{Proof.} Let $\bbX(\tau)\in q^\lambda\bbC^d[q^{-1},q]]$, i.e. $\bbX(\tau)=q^\lambda
\sum_{n=-N}^\infty\bbX_{(n)}q^n$ for some $N=N(\bbX)$, and let $\bbY(\tau)\in\cM^{-\lambda}_{2-w}(\rho^*)$,
and define a pairing $\langle \bbX,\bbY\rangle$ to be the $q^0$-coefficient $f_0$ of $\bbX(\tau)^t
\bbY(\tau)=\sum_{i=1}^d\bbX_i(\tau)\,\bbY_i(\tau)=\sum_{n=-N}^\infty f_nq^n$. Note that
$\langle \bbX,\bbY\rangle=\sum_{-N\le n\le 0}\bbX_{(n)}\bbY_{(-n)}$ depends only on the
coefficients $\bbX_{(n)}$ for $n\le 0$, since by hypothesis $\bbY_{(k)}=0$ for $k<0$.
In other words,  the pairing $\langle \bbX,\bbY\rangle$ depends only on $\bbY(\tau)$ and the
principal part $\cP_\lambda\bbX(q)$ of $\bbX(\tau)$.

If $\bbX(\tau)\in\cM^!_w(\rho)$, then $\bbX(\tau)^t\bbY(\tau)$ will lie in $\cM^!_2(1)$. Hence from the observation
after Proposition 3.4, in that case $\langle \bbX,\bbY\rangle$ will vanish. This means
the pairing $\langle \bbX,\bbY\rangle$ is a well-defined pairing between the cokernel of
$\cP_{\lambda;(\rho,w)}$, and $\cM^{-\lambda}_{2-w}(\rho^*)$.

Let $\bbY^1(\tau),\ldots,\bbY^m(\tau)$ be a basis of $\cM^{-\lambda}_{2-w}(\rho^*)$. We can require that
this basis be triangular in the sense that for each $1\le i\le m$ there is an $n_i,k_i$ so that
the coefficient  $\bbY^j_{(n_i)\,k_i}=\delta_{ji}$ for all $i,j$. Indeed, choose any $n_1,k_1$ such that $\bbY^1_{(n_1)\,k_1}\neq 0$, and rescale $\bbY^1(\tau)$ so that coefficient equals 1. Subtract if necessary a multiple of $\bbY^1(\tau)$ from
 the other $\bbY^i(\tau)$ so that $\bbY^i_{(n_1)\,k_1}=0$. Now, repeat: choose any $n_2,k_2$
 such that $\bbY^2_{(n_2)\,k_2}\ne 0$, etc. If we take take $\bbX^i(\tau)$ to be $q^{-n_i}e_{k_i}$
 (i.e. all coefficients vanish except one coefficient in one component), then $\langle\bbX^i,\bbY^j
 \rangle=\delta_{ij}$. This form of nondegeneracy means that dim$\,\cM^{-\lambda}_{2-w}(\rho^*)
 \le \mathrm{dim\,coker}\,\cP_{\lambda;(\rho,w)}$. 
 
 Repeating this argument in the dual direction, more precisely replacing $\rho,w,\lambda$ with $\rho^*,2-w,
 -1_d-\lambda$ respectively, gives us the dual inequality dim$\,\cM^{\lambda+1_d}_w(\rho)\le\mathrm{dim\,coker}\,\cP_{-1_d-\lambda;(\rho^*,2-w)}$. However, from Proposition 3.4 we know that  $c_{(\rho^*,2-w)}
 =\mathrm{Tr}\,(-1_d-\Lambda)=-d-c_{(\rho,w)}$, so the index formula \eqref{index} gives us both
 \begin{align}\mathrm{dim}\,\cM_w^{\lambda+1_d}(\rho)-\mathrm{dim\,coker}\,\cP_{\lambda;(\rho,w)}=&
 c_{(\rho,w)}-\mathrm{Tr}\,\lambda\,,\\
 \mathrm{dim}\,\cM_{2-w}^{-\lambda}(\rho^*)-\mathrm{dim\,coker}\,\cP_{-1_d-\lambda;(\rho^*,2-w)}=&
 -d-c_{(\rho,w)}+d+\mathrm{Tr}\,\lambda\,.\end{align}
Adding these gives 
$$ \mathrm{dim}\,\cM_w^{\lambda+1_d}(\rho)+\mathrm{dim}\,\cM_{2-w}^{-\lambda}(\rho^*)=\mathrm{dim\,coker}\,\cP_{\lambda;(\rho,w)}+\mathrm{dim\,coker}\,\cP_{-1_d-\lambda;(\rho^*,2-w)}\,.$$ Together with the two inequalities, this shows
 that the dimensions of $\cM^{-\lambda}_{2-w}(\rho^*)$ and coker$\,\cP_{\lambda;(\rho,w)}$ match.
 Hence  the  
pairing $\langle\bbX,\bbY\rangle$ is nondegenerate and establishes the isomorphism \eqref{eq:serre2}. Equation \eqref{eq:index2}
now follows immediately from the index formula \eqref{index}.
\qed\medskip%\quad\textbf{QED} to Theorem 4.3

We suggest calling part (a) \textit{Serre duality} because ker$\,\cP_\lambda$ has an
interpretation as $H^0(B_\Gamma;\cA)$ for some space $\cA=\cA_{\lambda;(\rho,w)}$ of meromorphic functions
on which $\Gamma$ acts by $(\rho,w)$, while coker$\,\cP_\lambda$, being the obstruction
to finding meromorphic sections of our $(\rho,w)$-vector bundle, should have an
interpretation as some $H^1$. The shifts by $1_d$ and 2 would be associated
to the canonical line bundle. This interpretation is at this point merely a heuristic, however.

 Compare \eqref{eq:serre2} to Theorem 3.1 in
 \cite{bib:Bor}. There, Borcherds restricts attention
to weight $w=k$ a half-integer, groups commensurable to SL$_{2}(\bbZ)$,
representations $\rho$ with finite image, and $\lambda=\lambda^{hol}$. 
The assumption of finite image was
essential to his proof. 
Our proof extends to arbitrary $T$ and arbitrary genus-0 groups \cite{bib:Ga} (so by inducing the representation,
it also applies to any finite-index subgroup of a genus-0 group). 

  The most important special case of the dimension 
formula \eqref{eq:index2} is $\lambda=\lambda^{hol}$,  which  relates the dimensions of
holomorphic vvmf with those of vector-valued cusp forms.  Compare \eqref{eq:index2} to  \cite{bib:Sk}, where Skoruppa  obtained the formula assuming $2w\in\bbZ$,
and that $\rho$ has finite image. His proof used the Eichler-Selberg trace
formula. Once again we see that these results hold in much greater generality.
Of course thanks to the induction trick, \eqref{eq:index2} also gives
the dimension formulas for spaces of e.g. holomorphic and cusp forms, for any 
finite index subgroup of $\bG$.

Incidentally, it is possible to have nonconstant modular functions, holomorphic everywhere in
$\uhp^*$, for some multipliers with infinite image (see e.g. Remark 
2 in Section 3 of  \cite{bib:KM}).    On the other hand,
Lemma 2.4 of \cite{bib:KM} prove that for \textit{unitary} multipliers $\rho$, there are no nonzero
 holomorphic vvmf $\FA(\tau)\in \mathcal{M}^{\lambda^{hol}}_{w}(\rho)$ of
 weight $w<0$.  This implies there are no vector-valued cusp forms of weight $w\le 0$, for
 unitary $\rho$.

\subsection{The $q$-expansion coefficients}

Modular forms --- vector-valued or otherwise --- are most important for their $q$-expansion
coefficients. In this subsection we study these. A special case of Proposition 3.5(b) (namely,
$w=0$ and $\bbQ_\bbX=\bbQ$) is given in \cite{bib:AM}, but our proof generalises without
change to $T$ nondiagonalisable and to any genus-0 group \cite{bib:Ga}.

As a quick remark, note from Theorem 3.3(b) that the coefficients $\bbX_{(n)\,i}$ of a vvmf $\bbX\in
\cM^!_w(\rho)$ will all lie
in the field generated over $\bbQ$ by the entries of $\Lambda_w$ and $\chi_w$, as well as
the coefficients $\bbX_{(n)\,i}$  of the principal part
$\cP_\Lambda(\bbX)$, where $\Lambda$ is bijective and as always $\chi=\Xi_{(1)}$.  

Let $\cF$ denote the span of all $p(\tau)\,q^u h(q)$, where $p(\tau)\in\bbC[\tau]$, $u\in\bbC$,
and $h\in \bbC[[q]]$ is holomorphic at $q=0$. The components 
of any vvmf (including `logarithmic', where $T$ is nondiagonalisable) must lie in $\cF$. 
$\cF$ is a ring, where terms simplify in the obvious way (thanks to the series $h$'s converging
absolutely). $\cF$ is closed under both differentiation and $\tau\mapsto\tau+1$. The key to analysing
functions in $\cF$ is the fact that they are equal only when it is obvious that they are equal: 

\medskip\noindent\textbf{Lemma 3.2.} \textit{Suppose $\sum_{i=1}^np_i(\tau)\,q^{u_i} h_i(q)=0$
for all  $q$ in some sufficiently small disc about $q=0$, 
where each $p_i(\tau)\in\bbC[\tau]$, $u_i\in\bbC$, $0\le \mathrm{Re}\,u_i<1$,
and $h_i\in \bbC[[q]]$ is holomorphic at $q=0$, and all $u_i$ are pairwise distinct. Then for each $i$, either $p_i(\tau)$ is identically 0
or $h_i(q)$ is identically 0 (i.e. all coefficients of $h_i(q)$ vanish).}

\medskip\noindent\textit{Proof.} For each $v\in\bbC$, define an operator $\cT_v$ on $\cF$ by
$(\cT_vg)(\tau)=g(\tau+1)-e^{2\pi\Ii v}g(\tau)$. Note that $p(\tau)q^uh(q)$ lies in the kernel of
$\cT_u^k$ if (and as we will see only if) the degree of the polynomial $p(\tau)$ is $<k$. 

Assume for contradiction that no $p_i(\tau)$ nor $h_i(q)$ are  identically 0. Let $k_i$ be 
the degree of $p_i(\tau)$. Apply $\cT^{k_d+1}_{u_d}\circ\cdots\circ\cT^{k_2+1}_{u_2}\circ
\cT^{k_1}_{u_1}$ to  $\sum_{i=1}^np_i(\tau)\,q^{u_i} h_i(q)=0$ to obtain $aq^{u_1}h_1(q)=0$
for some nonzero $a\in\bbC$ and all $q$ in that disc. This forces $h_1(q)\equiv 0$, a contradiction.
\qed

\medskip\noindent\textbf{Proposition 3.5.} \textit{Let $(\rho,w)$ be admissible and $T$
diagonal, and choose any vvmf $\bbX(\tau)\in\cM_w^!(\rho)$. Write $\bbX(\tau)=q^\lambda\sum_{n=0}^\infty\bbX_{(n)}q^n$.}

\smallskip\noindent\textbf{(a)} \textit{Let $\sigma$ be any field automorphism
of $\bbC$, and for each $1\le i\le d$ define $\bbX_i^\sigma(\tau)=q^{\sigma(\lambda_{ii})}
\sum_{n=0}^\infty\sigma(\bbX_{(n)\,i})\,q^n$. Then
$\bbX^\sigma(\tau)\in \cM_{\sigma w}^!(\rho^\sigma)$, where $(\rho^\sigma,\sigma w)$ is admissible,
$T^\sigma=e^{2\pi\Ii\,\sigma\lambda}$, $e^{\pi\Ii\sigma w/2}S^\sigma$ is conjugate to $e^{\pi\Ii w/2}
S$ and  $e^{2\pi\Ii\sigma w/3}U^\sigma$ is conjugate to $e^{2\pi\Ii w/3}U$. $\bbX^\sigma(\tau)$
is  $\sigma\lambda$-holomorphic iff $\bbX(\tau)$ is $\lambda$-holomorphic. Choose any
bijective exponent $\Lambda$ and fundamental matrix $\Xi(\tau)$ for $(\rho,w)$; then
$(\rho^\sigma,\sigma w)$
has bijective exponent $\sigma\Lambda$ and fundamental matrix $\Xi^\sigma(\tau)$, where 
$\sigma$ acts on $\Xi(\tau)$ column-wise.}

\smallskip\noindent\textbf{(b)} \textit{Let $\bbQ_\bbX$ be the field generated over
$\bbQ$ by all Fourier coefficients $\bbX_{(n)}$ of $\bbX(\tau)$. Assume the components $\bbX_i(\tau)$ of $\bbX(\tau)$
are linearly independent over $\bbC$. Then both the weight $w$ and all exponents $\lambda_i$
 lie in $\bbQ_\bbX$.}\medskip
 
 \noindent\textit{Proof.} Start with (a), and consider first $w=0$.  Write $\chi=\Xi_{(1)}$. Then ${\Xi}
 (\tau)$ obeys the differential equation \eqref{fundeq1}. Recall that $J,E_2,E_4,E_6,\Delta$ all have coefficients
 in $\bbZ$ and hence are fixed by $\sigma$. Then ${\Xi}^{\sigma}(\tau)$ formally satisfies
\eqref{fundeq1} with $\Lambda_0,\chi_0$ replaced with $\sigma(\Lambda_0),\sigma(\chi_0)$ (recall
that \eqref{fundeq1} is equivalent to the recursions \eqref{eq:recurs}).
Therefore ${\Xi}^{\sigma}(\tau)$ is a fundamental solution of that differential equation, with entries meromorphic
in $\uhp^*$. 

In fact, thanks to the Fuchsian equation \eqref{fundeq1}, the only possible poles of $\Xi^\sigma
(\tau)$ are at the cusps
or elliptic fixed points. The behaviour at the elliptic points is easiest to see from \eqref{fundeq2};
in particular, $S^\sigma$ and $U^\sigma$ are conjugate to $e^{2\pi\Ii\cA^\sigma_2}$  and 
$e^{2\pi\Ii\cA_3^\sigma}$ respectively. But $\cA_j^\sigma=\sigma\cA_j$ entry-wise and both $\cA_j$ are
diagonalisable with rational eigenvalues, so $\cA_j^\sigma$ is conjugate to $\cA_j$
(as it has the identical eigenvalue multiplicities). 
 
To generalise to arbitrary weight $w$, it's clear from Lemma 3.1 that $(\Delta^{w/12})^\sigma(\tau)
=\Delta^{\sigma w/12}(\tau)$. A weight-$w$ fundamental matrix is $\Delta^{w/12}(\tau)$ times a weight-0 one.

Finally, if $\Fa(\tau)\in\cM^!_w(\rho)$, then there exists a polynomial ${p}(J)\in\bbC^d[J]$ such that 
$\Fa(\tau)={\Xi}(\tau){p}(J(\tau))$, so $\Fa^{\sigma}(\tau)={\Xi}^{\sigma}(\tau){p}(J(\tau))^\sigma$. 
Since $\sigma$ fixes the coefficients of  $J(\tau)$, $\Fa^{\sigma}(\tau)$ manifestly lies in
$\cM_{\sigma w}^{!}(\rho^\sigma)$.

\smallskip
Now turn to (b). Suppose first that there is some entry $\Lambda_{ii}\not\in\bbQ_{\Fa}$.
Then there exists some
field automorphism $\sigma$ of $\bbC$ fixing $\bbQ_{\Fa}$
but with  $\delta:=\sigma \Lambda_{ii}
-\Lambda_{ii}$ nonzero  (see e.g. \cite{bib:AM} for a proof of why such a $\sigma$ exists).

From part (a),
$\Fa^{\sigma}(\tau)\in\cM_{\,\sigma w}^{!}(\rho^\sigma)$. Then $\Fa^{\sigma}(\tau)=q^{\sigma\Lambda-
\Lambda}\Fa(\tau)$. Choose any $\gamma=\left({a\atop c}{b\atop d'}\right)\in\bG$ with $c\ne 0$
(because $\bG$ is Fuchsian of the first kind it will have many such $\gamma$). Then $\tau\mapsto-1/\tau$ gives
\begin{equation}\label{eq:AMequ}
S^\sigma\Fa^{\sigma}(\tau)=\tau^{w-\sigma w}\exp(-2\pi\Ii\, (\sigma\Lambda-\Lambda)/\tau)\,S\,\Fa(\tau)\,.\end{equation}

Write $f(\tau),g(\tau)\in \cF$ for the $i$th entry of $S\Fa$ and of $S^\sigma\bbX^\sigma(\tau)$,
respectively. 
Because the entries of $\Fa(\tau)$ are linearly independent, $f(\tau)$ (and also $g(\tau)$) are nonzero.  Write $w'=w-\sigma w$. By 
\eqref{eq:AMequ}, $g(\tau)=\tau^{w'}\exp(-2\pi\Ii \delta/\tau)\,{g}(\tau)$.
  Like all entries of $S^\sigma\Fa^\sigma(\tau)$,
$g(\tau)$ is killed by the order-$d$ differential operator 
\begin{equation}L_{S^\sigma\Fa^\sigma}=\sum \widetilde{h}_l(\tau)\left(
\frac{\mathrm{d}}{2\pi\Ii
\mathrm{d}\tau}\right)^l\nonumber\end{equation}
 obtained from Lemma 2.2(b) by expanding out each $\Der_{\sigma w}^l$;
note that each $\widetilde{h}_l(\tau)\in q^v\bbC[[q]]$ for some $v\in\bbC$, being a combination
of the modular forms $h_l(\tau)$ of Lemma 2.2(b) and various derivatives of $E_2(\tau)$.
The product rule and induction on $l$ gives \begin{equation}\label{prodrule}
 \left(
\frac{\mathrm{d}}{2\pi\Ii
\mathrm{d}\tau}\right)^lg(\tau)=  \tau^{w'}\exp(-2\pi\Ii\, \delta/\tau)\sum_{k=0}^l\frac{p_{l,k}(\tau)}{\tau^{2l-2k}} \left(
\frac{\mathrm{d}}{2\pi\Ii
\mathrm{d}\tau}\right)^k{f}(\tau)\,,\end{equation}
 where $p_{l,k}(\tau)$ is a polynomial in $\tau$ of degree $\le l-k$ and $p_{l,0}(\tau)$ has
  nonzero constant term $(2\pi\Ii\delta)^{l-k}$.  Multiplying $L_{S^\sigma\Fa^\sigma}g=0$ by $\tau^{2d}\tau^{-w'}\exp(2\pi\Ii \delta/\tau)$, we obtain
  \begin{equation}\sum_{l=0}^d\widetilde{h}_l(\tau)\sum_{k=0}^l\tau^{2k}p_{l,k}(\tau) 
  \left(
\frac{\mathrm{d}}{2\pi\Ii
\mathrm{d}\tau}\right)^k{f}(\tau)=0\,.\label{eq:mess}\end{equation}
Now, the derivatives $ \left(\frac{\mathrm{d}}{2\pi\Ii\mathrm{d}\tau}\right)^k{f}(\tau)$ are 
manifestly $q$-series (i.e. functions in $\cF$ where the polynomial parts $p_i(\tau)$ are all
constant), so we see from Lemma 3.2 that (regarding \eqref{eq:mess} as a polynomial in
$\tau$ with $q$-series coefficients) the $\tau^0$ coefficient of \eqref{eq:mess} must itself
vanish. This is simply $\widetilde{h}_d(\tau)\,(2\pi\Ii\delta)^df(\tau)=0$, where $\widetilde{h}_d
(\tau)=\mathrm{Wr}(\bbX^\sigma)(\tau)\ne 0$ by Lemma 2.2(a). This forces $\delta=0$,
i.e. $\sigma\Lambda_{ii}=\Lambda_{ii}$. 
 
Therefore all entries $\Lambda_{ii}$ must lie in $\bbQ_{\Fa}$.
Suppose now that the weight $w$  doesn't lie in 
$\bbQ_{\Fa}$,
and as above choose a 
field automorphism $\sigma$ of $\bbC$ fixing $\bbQ_{\Fa}$ but with $\sigma
w-w\not\in\bbZ$  (if $w$ lies in an algebraic extension of $\bbQ_{\Fa}$ then this is 
automatic, while if $w$ lies in a transcendental extension we
can select $\sigma w$ to likewise be an arbitrary transcendental). The remainder of the argument is as above: $\frac{p_{l,k}(\tau)}{\tau^{2l-2k}}$
in \eqref{prodrule} is now replaced with $\frac{c_{l,k}}{\tau^{l-k}}$ for some $c_{l,k}
\in\bbC$, and $c_{l,0}=w'(w'-1)\cdots(w'-l+1)\ne 0$.  Multiplying $L_{S^\sigma\Fa^\sigma}g=0$ by 
$\tau^{d}\tau^{-w'}$, we see that $\widetilde{h}_d(\tau){c_{d,0}}{f}(\tau)=0$,
likewise impossible unless $w\in\bbQ_\bbX$.
\qed%   \qquad \textbf{QED} to Corollary 3.8
   \medskip

The obvious Galois action on representations, namely $(\sigma\rho)(\gamma)_{ij}=\sigma(\rho
 (\gamma)_{ij})$, is unrelated to this $\rho^\sigma$. It would be interesting though to understand 
 the relation between the vvmf of  $\sigma\rho$ and those of $\rho$. 
 
For vvmf with rank-1 multipliers, we have $T^\sigma=e^{2\pi\Ii\sigma\lambda}$, 
$e^{\pi\Ii\sigma w/2}S^\sigma=e^{\pi\Ii w/2}S$, and
 $e^{2\pi\Ii \sigma w/3}U^\sigma=e^{2\pi\Ii w/3}U$. When the rank $d$ is greater than 1, however, the precise formula for $S^\sigma$ 
 and $U^\sigma$ is delicate. For example, we learn in Section 4.2 that when $d=2$,
 the explicit relation between $S^\sigma$ and $S$ involves the relation between the Gamma
function values $\Gamma(\sigma \Lambda_{11})$ and $\Gamma(\Lambda_{11})$. 
 
The next result is formulated in terms of certain modules $\cK$. One important example
is $\cK=\bbK[q^{-1},q]]$ for any subfield $\bbK$ of $\bbC$. 
Another example  is the subset $\cK$ of $f\in\bbQ[[q]]$
with bounded denominator, i.e. for which there is an $N\in\bbZ_{>0}$  such that
$Nf\in\bbZ[[q]]$. Both examples satisfy all conditions of Proposition 3.6.
This latter example can be refined in several ways, e.g. by fixing from the start a set $P$ of
primes and requiring that the powers of the primes in $P$ appearing in the denominators be bounded, but primes $p\not\in P$
be unconstrained. 
 
\medskip
\noindent\textbf{Proposition 3.6.} 
\textit{Suppose $(\rho,w)$ is admissible, $\rho$ is irreducible and $T$ is diagonal, and choose any exponent $\lambda$. Let $\bbK$ be any subfield of $\bbC$ and let $\mathcal{K}\subset\bbK[q^{-1},q]]$
be any module over both $\bbK$ and $\bbZ[q^{-1},q]]$  (where both of these  act by multiplication), such that $\frac{\mathrm{d}}{\mathrm{d}q}\cK\subseteq
\cK$. Let $\mathcal{M}_w^{!\ \mathcal{K}}(\rho)$ denote the intersection
$\mathcal{M}^!_w(\rho)\cap q^\lambda\mathcal{K}$. Then
the following are equivalent:}
\begin{itemize}
\item[(i)] $\ \mathcal{M}_w^{!\ \mathbb{K}}(\rho)\ne 0$;

\item[(ii)] \textit{\ $\mathcal{M}_{w+2k}^{!\ \mathbb{K}}(\rho)\ne 0$ for all $k\in\bbZ$};

\item[(iii)] $\ \mathrm{span}_\bbC\,\mathcal{M}_w^{!\  \mathbb{K}}(\rho)=\mathcal{M}^!_w(\rho)$;

\item[(iv)] \textit{\ for any bijective exponent $\Lambda$, all entries $\Xi_{ij}(\tau)$ of the corresponding 
fundamental matrix $\Xi(\tau)$ lie in $q^{\Lambda_{ii}}\cK$;}

\item[(v)] \textit{\ for any exponent $\lambda$, the space $\mathrm{span}_\bbC\left(\cM^{\lambda}(\rho)\cap q^\lambda
\mathcal{K}\right)=\cM^{\lambda}(\rho)$.}
\end{itemize}
\medskip

\noindent\textbf{Proof.} Assume (i) holds, and choose any nonzero $\bbX(\tau)\in  \mathcal{M}_w^{!\ \cK}(\rho)$. Then $\Delta^hE_4^iE_6^j\bbX\in q^\lambda\cK$ for any $h\in\bbZ$ and $i,j\in\bbZ_{\ge 0}$, which gives us (ii). 
That $\rho$ is irreducible forces the components of $\bbX(\tau)$ to be
linearly independent over $\bbC$.
Then Proposition 3.5 would require all $w,\lambda_{ii}\in\bbK$, and hence $q^\lambda\cK^d$ 
is mapped into itself by $\DO_{w',i}$ for any $w'\in w+2\bbZ$. Moreover, Proposition 3.2 applies and $\cM_w^!(\rho)=
\bbC[J,\DO_1,\DO_2,\DO_3]\bbX(\tau)$, which gives us (iii). Let $\Lambda$ be any bijective exponent. 
Then each column of $\Xi(\tau)$ is a linear combination over $\bbC$ of finitely many
vvmf in $\cM^{!\ \cK}_w(\rho)$; but all of these vectors are uniquely determined by their principal
parts $\cP_\Lambda$, which are all in $\bbK^d[q^{-1}]$, so that linear combination must have a 
solution over the field $\bbK$. This gives us (iv).
To get (v), note that in the proof of Theorem 3.4(a) we may choose our $\bbX^{(i)}(\tau)$ to lie in
$q^\lambda\cK$, since all we require of them is a linear independence condition.
\qed%\textbf{QED} to Theorem 5.2
\medskip

This proposition is relevant to
the study of 
 modular forms for noncongruence subgroups of $\Gamma$. 
A conjecture attributed to Atkin--Swinnerton-Dyer \cite{bib:AS} states that a (scalar)
modular form for some subgroup of $\Gamma$ will have bounded denominator only
if it is a modular form for some congruence subgroup. More generally, it is expected
that a vvmf $\bbX(\tau)$ for $\Gamma$, with entries $\bbX_i(\tau)$ linearly independent over $\bbC$ and
 with coefficients $\bbX_{(n)\,i}$ all in $\bbQ$, will have bounded denominators
only if its weight $w$ lies in $\frac{1}{2}\bbZ$ and the kernel of $\varmu_{-w}\otimes\rho$
is a congruence subgroup. If the kernel is of infinite index, then it is expected that infinitely
many distinct primes will appear in denominators of coefficients.

\subsection{Semi-direct sums and exactness}

If $\rho$ is a direct sum $\rho'\oplus\rho''$, then trivially its vvmf are comprised of those of
$\rho'$ and $\rho''$. But what if $\rho$ is a semi-direct sum?
 
Consider a short exact sequence
\begin{equation}\label{eq:ses} 0\rightarrow U\stackrel{\iota}{\rightarrow} V\stackrel{\pi}{\rightarrow} W\rightarrow 0\end{equation}
of finite-dimensional $\overline{\bG}$-modules. We can consider an action at arbitrary weight by
tensoring with $\varmu_w$. Choose a basis for $V$ in which
$\rho_V=\left({\rho_U\atop 0}{*\atop \rho_W}\right)$. In terms of this basis, $\bbX_V\in \cM_{ w}^{!}
(\rho_V)$ implies $\bbX_V=\left({\bbX_U\atop \bbX_W}\right)$ where $\bbX_W\in\cM_{w}^{!}(\rho_W)$,
and $\left({\bbX_U\atop 0}\right)\in\cM_{w}^{!}(\rho_V)$ iff $\bbX_U\in\cM_{w}^{!}
(\rho_U)$.  For this basis we have the natural embedding
 $\iota'(\bbX_U)= \left({\bbX_U\atop 0}\right)$ and projection $\pi'\left({\bbX_U\atop \bbX_W}\right)=\bbX_W$. 
We write \eqref{eq:ses} as $\rho_V=\rho_U\sdsum\rho_W$.

The fundamental objects in the theory of vvmf are the functors $(\rho,w)\mapsto\cM^!_w(\rho)$,
$(\rho,w;\lambda)\mapsto \cM^{\lambda}_w(\rho)$ and $(\rho,\lambda)\mapsto \cM^{\lambda}(\rho)$, attaching spaces of vvmf to
multipliers, weights and exponents. 
Marks--Mason \cite{bib:MM}  suggest considering the effects of these functors on
\eqref{eq:ses}. It is elementary that they are left-exact, i.e. that  $\rho_V=\rho_U\sdsum\rho_W$
trivially implies 
\begin{equation}\label{eq:left}
0\rightarrow \cM_{w}^{\lambda_U}(\rho_U)\stackrel{\iota'}{\rightarrow}\cM_{w}^{\lambda_V}(\rho_V)
\stackrel{\pi'}{\rightarrow}
\cM_{w}^{\lambda_W}(\rho_W)\end{equation}
and similarly for $\cM_w^!$ etc. However \cite{bib:MM} found 1-dimensional $U$ and $W$
such that the functor $\rho\mapsto\cM^{\lambda^{hol}}(\rho)$ is not right-exact.

Thanks to Theorem
3.5, we can generalise and quantify this discrepancy. 
We thank Geoff Mason for suggesting the naturalness of explaining the failure of
right-exactness with a long exact sequence.

\medskip\noindent\textbf{Theorem 3.6.} \textit{Write $\rho_V=\rho_U\sdsum\rho_W$ as in \eqref{eq:ses}, and suppose $(\rho_V,w)$ is admissible and $T_V$
diagonal. Choose any
 exponent $\lambda_V=\mathrm{diag}(\lambda_U,
\lambda_W)$. Then we obtain the exact sequences of $\bbC$-spaces} 
\begin{eqnarray}\label{eq:ses2}
0&\rightarrow& \cM_{w}^{!}(\rho_U)\stackrel{\iota'}{\rightarrow}\cM_{w}^{!}(\rho_V)
\stackrel{\pi'}{\rightarrow}
\cM_{w}^{!}(\rho_W)\rightarrow 0\,,\\
0\rightarrow \cM_U\stackrel{\iota'}{\rightarrow}\cM_V&\stackrel{\pi'}{\rightarrow}&
\cM_W\stackrel{\delta}{\rightarrow}\mathrm{coker}\,\cP_U%\nonumber\\ &&\qquad
\cong (\cM_{U\ast})^*\stackrel{\pi''*}{\rightarrow}(\cM_{V\star})^*
\stackrel{\iota''*}{\rightarrow}(\cM_{W\ast})^*\rightarrow 0\label{eq:les}
\end{eqnarray}
\textit{
where we write $\cP_U=\cP_{\lambda_U-1_{d_U};(\rho_U,w)}$, $\cM_U=\cM^{\lambda_U}_w
(\rho_U)$, $\cM_{U\ast}=\cM_{2-w}^{1_{d_U}-\lambda_U}(\rho_U^*)$ etc. $\iota''{}^*$ and $\pi''{}^*$ are restrictions of the dual (transpose)
maps. The isomorphism in \eqref{eq:les} is \eqref{eq:serre2}, and the connecting map $\delta$ is defined in the proof. Moreover, for any bijective exponents $\Lambda_U$ of $(\rho_U,w)$ and 
$\Lambda_W$ of $(\rho_W,w)$, diag$(\Lambda_U,\Lambda_W)$ is bijective for $(\rho_V,w)$.}

\medskip\noindent\textit{Proof.} First let's prove \eqref{eq:ses2}. Let $\cM^!_V,\cM^!_U,\cM^!_W$ denote $\cM_{w}^{!}(\rho_V),
\cM_{w}^{!}(\rho_U),\cM_{w}^{!}(\rho_W)$ respectively. Choose any bijective exponents
$\Lambda_U,\Lambda_W$ and define  $\Lambda_V=\mathrm{diag}(\Lambda_U,
\Lambda_W)$, $\cP'_V,\cP'_U,\cP'_W$ for $\cP_{\Lambda_V}$ etc. If $\left({\bbX_U\atop \bbX_W}\right)\in\mathrm{ker}\,
\cP'_V$, then $\bbX_W(\tau)\in\mathrm{ker}\,\cP'_{W}$ and hence
$\bbX_W(\tau)=0$ because $\Lambda_W$ is bijective. This means $\bbX_U(\tau)\in\cM^!_U$, so also
$\bbX_U(\tau)\in\mathrm{ker}\,\cP'_{U}$ and $\bbX_U(\tau)=0$ because $\Lambda_U$ is bijective.
Therefore $\cP'_{V}$ is injective. From the first line of \eqref{cvalue} we know $c_{(\rho_V,w)}=
c_{(\rho_U,w)}+c_{(\rho_W,w)}$ (since Tr$\,S_V=\mathrm{Tr}\,S_U+\mathrm{Tr}\,S_W$ etc); the 
index formula \eqref{index} then implies coker$\,\cP'_{V}=0$ and hence $\Lambda_V$ is bijective. 
In order to establish \eqref{eq:ses2}, 
only the surjectivity of $\pi'$ needs to be shown, but this follows from the surjectivity of $\cP'_V$
together with the injectivity of $\cP'_{W}$.

Now let's turn to  \eqref{eq:les}. Most of this exactness again comes from \eqref{eq:left}: the dual of
$0\rightarrow\cM_{W*}\rightarrow\cM_{V*}\rightarrow\cM_{U*}$ gives the second half of
\eqref{eq:les}. The connecting map $\delta$ is defined as follows.  Given any 
$\bbX_W\in\cM_W$, exactness of \eqref{eq:ses2} says that there is an $\bbX_U$, unique mod
$\cM^!_U$, such that $\left({\bbX_U\atop \bbX_W}\right)\in\cM^!_V$. The connecting map $\delta:\cM_{W}\rightarrow\mathrm{coker}\,\cP_U$ sends $\bbX_W$ to $\cP_U(\bbX_U)+\mathrm{Im}\,\cP_U$.
Exactness at $\cM_W$ is now clear: if $\bbX_W\in\mathrm{ker}\,\delta$ then $\cP_U(\bbX_U)
=\cP_U(\bbX'_U)$ for some $\bbX'_U\in\cM_U$, so $\left({\bbX'_U\atop 0}\right)$ and hence
$\left({\bbX_U-\bbX'_U\atop \bbX_W}\right)$ both lie in $\cM^!_V$. However the latter
manifestly lies in the kernel of $\cP_V$, so $\bbX_W=\pi'\left({\bbX_U-\bbX'_U\atop \bbX_W}\right)$ as desired.
To see exactness at $\pi''{}^*$, suppose $f_U\in\mathrm{ker}\,\pi''{}^*$. Then $f_U$ is a functional
on $\cM_{U*}$ and the associated functional $\left({\bbY_U\atop \bbY_W}\right)\mapsto f_U(\bbY_U)$
on $\cM_{V*}$ is 0. Thanks to \eqref{eq:serre2}, these functionals can be expressed as 
$\bbY_U\mapsto\langle\bbP_U',\bbY_U\rangle$ and $\left({\bbY_U\atop \bbY_W}\right)
\mapsto \langle\left({\bbP_U\atop \bbP_W}\right),\left({\bbY_U\atop \bbY_W}\right)\rangle$ for
some $\bbP'_U(\tau)\in\bbC^{d_U}[q^{-1}]$ and $\left({\bbP_U\atop \bbP_W}\right)\in\bbC^{d_U+d_W}
[q^{-1}]$. The independence of the functional on
$\bbY_W$ means $\bbP_W(\tau)=0$, so we can take $\bbP_U(\tau)=\bbP_U'(\tau)$. That the functional on $\cM_{V*}$ is 0 means (again from \eqref{eq:serre2}) that
$\left({\bbP_U\atop 0}\right)=\cP_V\left({\bbX_U\atop \bbX_W}\right)$ for some
$\left({\bbX_U\atop \bbX_W}\right)\in\cM^!_V$. Moreover, $\bbX_W(\tau)\in\cM_W$, so $\bbP_U=
\delta(\bbX_W)$, as desired.\qed\medskip% \textbf{QED} to Theorem 4.4

Theorem 3.6 allows us to classify all bijective $\Lambda_V$. These are given by all exponents
diag$\,(\lambda_U,\lambda_W)$ such that $\cP_{\lambda_U}$ is injective, $\cP_{\lambda_W}$  
is surjective, and the connecting map $\delta:\mathrm{ker}\,\cP_{\lambda_W}\rightarrow
\mathrm{coker}\,\cP_{\lambda_U}$ is an isomorphism.

We can now quantify the failure of $\cM^{\lambda}(\rho)$ to be exact. 
 For each fixed $w$, the discrepancy is
\begin{equation}\mathrm{dim}\,\cM_U+\mathrm{dim}\,\cM_W-\mathrm{dim}\,\cM_V
=\mathrm{dim\,Im}\,\delta\,.\label{eq:discr}\end{equation}
For $w\ll 0$, $\cM^{\lambda}_w(\rho_W)=0$ so $\delta=0$ and the discrepancy is 0, while
for $w\gg 2$, then $\cM^{1_{d_U}-\lambda^{hol}}_{2-w}(\rho_U^*)=0$ so again
$\delta=0$ and the discrepancy is 0. Thus the total discrepancy, summed over all $w$, is
finite.

Let us recover in our picture  the calculation in Theorem 4 of \cite{bib:MM}.  Take $\lambda=\lambda^{hol}$. 
Consider  $\rho_V$ of the form $\varmu_{2a}\sdsum\varmu_{2b}$, where as always $\varmu_j$ has $T=
e^{2\pi\Ii j/12}$. Then Theorem 4 of \cite{bib:MM} says $\rho_V$ can be  indecomposable iff
$|a-b|=1$. As above, if $w<0$ then $\cM_W=0$ while if $w\ge 2$ then $\cM_{U*}=0$,
so only at $w=0$ can $\delta\ne 0$. We find that $w=0$ and $\cM_W,\cM_{U*}\ne0$ forces 
$b=0$ and $a=5$, in which case $\delta:\bbC\rightarrow\bbC$. Now, a bijective exponent for
$(\rho_V,0)$ is diag$\,(-\frac{7}{6},0)$ by Theorem 3.6 (and Section 4.1 below), so there 
must be a $\bbX_U(\tau)\in q^{-1/6}\sum_{n=0}^\infty c_nq^n$ such that $\left({\bbX_U\atop 1}
\right)\in\cM_V^!$, by
surjectivity. Indeed, $\delta(1)=\cP_U(\bbX_U)=c_0$. If $c_0=0$ then the Wronskian of 
$\left({c_1q^{5/6}+\cdots\atop 1}\right)$ would be a nonzero holomorphic modular form $cq^{5/6}+\cdots$
of weight 2 (the Wronskian is nonzero because $\rho_V$ is indecomposable). This is impossible (e.g. $\eta^{-20}$ times it would also be holomorphic but with trivial
multiplier and weight $-8$). Therefore $c_0\ne 0$, so $\delta\ne 0$ and the total
discrepancy is 1-dimensional.

\section{Effectiveness of the theory}

Explicit computations within our theory are completely feasible.
Recall from Theorem 3.3 that we have complete and explicit knowledge of the space $\cM_{w}^!({\rho})$
of weakly holomorphic vvmf, if we know the diagonal matrix $\Lambda$ and
the complex matrix $\mathcal{X}$. We know Tr$\,\Lambda$, and generically any matrix
with the right trace and with $e^{2\pi\Ii\Lambda}=T$ is a bijective exponent.
This $\mathcal{X}$ can be obtained in principle from $\rho$
and $\FB$ using the Rademacher expansion (see e.g. \cite{bib:DMMV}).
 However this series expansion for $\mathcal{X}$ converges notoriously slowly, and
obscures any properties of $\mathcal{X}$ that may be present (e.g.  integrality). Hence
other methods are needed for identifying $\mathcal{X}$. In this section we provide 
several examples, illustrating some of the ideas available. See 
\cite{bib:BG,bib:BG2} for further techniques and examples.

\subsection{One dimension}

It is trivial to solve the $d=1$ case \cite{bib:BG2}. 
Here, $\rho=\varmu_{u}$, and $u\in\bbC$, $0\le 
\mathrm{Re}\,{u}<1$, and the weight $w$ is required to be $w\in 12u+2\bbZ$. Write $w=12u+2j+12n$
where $0\le j\le 5$ and $n\in\bbZ$; then for 
$j= 0,1,2,3,4,5$ resp., the fundamental matrix $\Xi(\tau)$ for $\cM^!_{12u+2j+12n}({\varmu_{u}})$
is (recall Lemma 3.1):
\begin{align}\nonumber &\Delta^{u+n}(\tau)\,,\qquad\Delta^{u+n-5/6}(\tau)\,E_{10}(\tau)\,,\qquad\Delta^{u+n-2/3}(\tau)\,E_8(\tau)\,,\\
&\Delta^{u+n-1/2}(\tau)\,E_6(\tau)\,,\qquad\Delta^{u+n-1/3}(\tau)\,E_4(\tau)\,,\qquad\Delta^{u+n-7/6}(\tau)\,E_{14}(\tau)\,,
\nonumber\end{align}
 respectively, where $E_8=E_4^2$, $E_{10}=E_4E_6$ and $E_{14}=E_4^2E_6$.
The unique free generator of $\cM^{\lambda^{hol}}(\varmu_u)$ is $\Delta^{u}(\tau)$.
This means that dim$\,\cM_w^{\lambda^{hol}}(\varmu_u)$ equals the dimension of
the weight $w-12u$ subspace of $\frakm$, assuming as above that $0\le\mathrm{Re}\,u<1$.

\subsection{Two dimensions}

Much more interesting is $d=2$ (see e.g. \cite{bib:TW,bib:M2}). Let's start with weakly
holomorphic. As usual, it suffices to consider  weight $w=0$. We continue to require
that $T$ be diagonal, although we give a `logarithmic' example shortly.

The moduli space of equivalence classes of 2-dimensional representations
of $\overline{\Gamma}$ consists of 15 isolated points, together with 3 half-planes.
Each half-plane has 4 singularities: 2 conical singularities and 2 triple points. The
equivalence classes of irreducible representations with $T$ diagonalisable correspond 
bijectively with the regular
points on the 3 half-planes. The 15 isolated points are all direct sums $\rho_1\oplus\rho_2$
of 1-dimensional representations which violate the inequalities \eqref{eq:inequal}, and are
handled by Section 4.1. Each conical singularity corresponds to
an irreducible `logarithmic' representation with $T$ nondiagonalisable. Each triple point consists of a direct
sum $\rho\oplus\rho'$ as well as the semi-direct sums $\rho\sdsum\rho'$ and $\rho'\sdsum\rho$. 
In all cases except the 6 triple
points, the set of eigenvalues of $T$ uniquely determine the representation.
The 3 half-planes correspond to the 3 different choices of $(\alpha_i,\beta_j)$
possible at $d=2$ (recall \eqref{eq:inequal}).

More explicitly, one half-plane  has det$\,T=\xi_6=:\zeta$:
\begin{equation}
T=\left(\begin{array}{cc}z&0\\ 0&\zeta z^{-1}\end{array}\right)\ ,\qquad S=\left(
\begin{array}{cc}
\frac{\bar{\zeta}z}{z^2-\zeta}&y\\ \frac{(z^2-1)(z^2-\zeta^2)}{y
(z^2-\zeta)^2}&-\frac{\bar{\zeta}z}{z^2-\zeta}\end{array}\right)\,,
\end{equation}
for arbitrary $y,z\in\bbC$ provided $y\ne 0$ and
$z\not\in\{ 0,\pm \xi_{12}\}$. This is irreducible iff $z\not\in\{\pm 
1,\pm\zeta\}$. The redundant parameter $y$ is introduced for later convenience.
\textit{Irreducible} $\rho,\rho'$ with $z,z'$ related by 
$(zz')^2=\zeta$ are naturally isomorphic. Here, $\alpha_1=\beta_0=\beta_1=1$,
so Tr$\,\FB=-\frac{5}{6}$ and $\FB=\textrm{diag}(t,-\frac{5}{6}-t)$
where $z=e^{2\pi \Ii t}$ for some $t\not\equiv \frac{1}{12},\frac{7}{12}$ (mod 1). 
Then Theorem 3.3(c) tells us
\begin{equation}
\mathcal{X}=\left(\begin{array}{cc} 24\frac{t(60t-11)}{12t+5}&10368x\frac{t(2t+1)(3t+1)(6t+5)}{(12t+11)(12t+5)^2}\\
\frac{10368}{x(12t-1)}&-4\frac{(6t+5)(60t+61)}{12t+5}\end{array}\right)\,,
\end{equation}
for some $x\ne 0$ to be determined shortly. Equation \eqref{fundeq2}
is ideally suited to relate $x$ and $y$,
since at $d=2$ it reduces to the classical hypergeometric equation. We read off
from it the fundamental matrix
\begin{equation}\label{eq:fund2a}
%\begin{align}\label{eq:fund2a}
{\Xi}(\FC(\tau))= %\,&=\\ 
\left(\begin{array}{cc} f(t;\frac{5}{6};\FC(\tau))&
\mathcal{X}_{12}f(t+1;\frac{5}{6};\FC(\tau))
\\ \mathcal{X}_{21}f(\frac{1}{6}-t;\frac{5}{6};\FC(\tau))&
f(-\frac{5}{6}-t;\frac{5}{6};\FC)\end{array}\right)
\end{equation}
for $\FC(\tau)=J(\tau)/1728$, where we write 
\begin{equation}\label{eq:hypser}
f(a;c;\FC)=(-1728\FC)^{-a}F(a,a+\frac{1}{2};2a+c;\FC^{-1}) \end{equation}
for $F(a,b;c;z)=1+\frac{ab}{c}z+\cdots$ the hypergeometric series.
%\begin{equation}\label{eq:hypser}
%F(a,b;c;z)= \sum_{n=0}^\infty \frac{(a)_n(b)_n}{n!(c)_n}z^n  
%\end{equation}
Substituting $\FC(\tau)=J(\tau)/1728$ directly into \eqref{eq:hypser} and \eqref{eq:fund2a} gives the
$q$-expansion of ${\Xi}(\tau)$. The parameters $x,y$ appearing in $S$ and 
$\mathcal{X}$ can be related by the standard analytic continuation of 
$F(a,b;c;z)$ from $z\approx 0$ to $z\approx \infty$, which implies
\begin{equation}
f(a;c;\FC)=(1728)^{-a}\left\{\frac{\sqrt{\pi}\Gamma(2a+c)}{\Gamma(a+\frac{1}{2})\,
\Gamma(a+c)}-\frac{2\sqrt{\pi}\Gamma(2a+c)}{\Gamma(a)\,\Gamma(a+c-\frac{1}{2})}
\FC^{1/2}+\cdots\right\}\end{equation}
for small $|\FC|$. Hence
\begin{equation}
y=\frac{\sqrt{3}x}{1728}\frac{2^{2/3}}{432^{2t}}\frac{\Gamma(2t+\frac{5}{6})^2}{\Gamma(2t)\,\Gamma(2t+
\frac{2}{3})}\,.
\label{eq:xyclass1}\end{equation}
In particular, we see that when $\rho$ is irreducible, $\FB$ is bijective for a $\bbZ$
worth of $t$'s (i.e. the necessary conditions $e^{2\pi\Ii\FB}=T$ and Tr$\,\FB=-\frac{5}{6}$ are also
sufficient) --- at the end of Section 3.4 we call such $\rho$ \textit{tight}. This fact is generalised in Theorem 4.1 below.
There are  4 indecomposable but reducible $\rho$ here:

\smallskip \noindent $z=1$: this $\rho$ is the semi-direct sum 
$1\sdsum\varmu_{2}$ of 1 (the subrepresentation)
with  $\varmu_{2}$ (the quotient); in this case $\FB=\textrm{diag}(t,-\frac{5}{6}-t)$ is bijective
iff $t\in 1+\bbZ_{\ge 0}$;

\smallskip\noindent $z=-1$: here, $\rho=\varmu_{6}\sdsum\varmu_{8}$; its $\FB$ is bijective iff 
$t\in\frac{1}{2}+\bbZ_{\ge 0}$;

\smallskip\noindent $z=\zeta$: here, $\rho=\varmu_{2}\sdsum 1$; its $\FB$ is bijective iff $t\in
\frac{1}{6}+\bbZ_{\ge 0}$;

\smallskip\noindent $z=-\zeta$: here, $\rho=\varmu_{8}\sdsum\varmu_{6}$; its $\FB$ is 
bijective iff $t\in\frac{2}{3}+\bbZ_{\ge 0}$.
\smallskip

\noindent These restrictions on $t$ are needed to avoid the Gamma function poles in
\eqref{eq:xyclass1}. Nevertheless, the `missing' values of $t$ (apart from the forbidden
$t\equiv \frac{1}{12}$ (mod $\frac{1}{2}$)) are all accounted for: for example $t=0$ describes
the direct sum $1\oplus\varmu_{2}$ (the limiting case where $y\rightarrow 0$  
slowly compared with $t\rightarrow 0$),
while $t\in\bbZ_{< 0}$ recovers the $z^2=\zeta$ solution given above.

The free generators over $\frakm$ of $\cM^{\lambda^{hol}}(\rho)$ for these $\rho$ are now
easy to find. Fix $0\le\mathrm{Re}\,t<1$, and note that $\lambda^{hol}=\mathrm{diag}
(t,\frac{1}{6}-t)$ when Re$\,t\le\frac{1}{6}$, and otherwise $\lambda^{hol}=\mathrm{diag}
(t,\frac{7}{6}-t)$ .Consider first the case where $\rho$ is irreducible; then $\rho$ is
tight and Proposition 3.3(b) tells us dim$\,\cM_w^{\lambda^{hol}}(\rho)$ for all $w<0$.
In particular, if Re$\,t\le \frac{1}{6}$ then $w^{(1)}=0$ and $w^{(2)}=2$, and $\bbX^{(1)}(\tau)$
is the first column of $\Xi(\tau)$ given above; if instead Re$\,t>\frac{1}{6}$
then  $w^{(i)}=6,8$, $\bbX^{(1)}(\tau)$ is the first column of $\Xi(\tau)$ at $w=6$ for $\Lambda=\mathrm{diag}
(t,\frac{1}{6}-t)$. In both cases, $\bbX^{(2)}(\tau)=\Der\bbX^{(1)}(\tau)$.

The holomorphic vvmf for the 2-dimensional indecomposable representations was discussed 
at the end of Section 3.7, and we find that for our four such $\rho$, dim$\,\cM^{\lambda^{hol}}_w(\varmu_{2a}\sdsum\varmu_{2b})=
\mathrm{dim}\,\cM^{\lambda^{hol}}_w(\varmu_{2a})+\mathrm{dim}\,\cM^{\lambda^{hol}}_w(\varmu_{2b})$. The 1-dimensional case was worked out in Section 4.1, and we find that
in all cases $\{w^{(1)},w^{(2)}\}=\{2a,2b\}$ with $\bbX^{(i)}(\tau)$ given by the appropriate column of 
the fundamental matrix at $w=2a$ and $2b$ respectively.

The choice $z=\pm\zeta$, i.e. $t\equiv \frac{1}{12}$ (mod $\frac{1}{2}$), corresponds here to two logarithmic representations.
Consider for concreteness $z=\xi_{12}$. A weakly holomorphic vvmf for it is $\eta(\tau)^2
\left({\tau\atop 1}\right)$. This generates all of $\cM^!_0$, using $\bbC[J]$ and the differential operators $\nabla_i$; together with 
$$q^{1/12}\left(\begin{matrix}\pi\Ii q^{-1}-242\pi\Ii+(-140965\pi\Ii-55440\tau)q+\cdots\\ -55440q+\cdots
\end{matrix}\right)$$
it freely generates $\cM^!_0$ over $\bbC[J]$. The 5 other 2-dimensional logarithmic representations
correspond to this one tensored with a $\overline{\bG}$ character.  The free basis for
holomorphic vvmf is $\eta(\tau)^2
\left({\tau\atop 1}\right)$ and its derivative.

We can see the Galois action of Section 3.6 explicitly here: $\sigma$ takes 
$t\mapsto\sigma t$ and $x\mapsto \sigma x$, and it keeps one inside this connected
component.

\smallskip Another class of two-dimensional representations has det$\,T=-1$:
\begin{equation}
T=\left(\begin{array}{cc}z&0\\ 0&- z^{-1}\end{array}\right)\ ,\qquad S=\left(
\begin{array}{cc}
\frac{-z}{z^2+1}&y\\ \frac{z^4+z^2+1}{y
(z^2+1)^2}&\frac{z}{z^2+1}\end{array}\right)\,,
\end{equation}
for arbitrary $y,z\in\bbC$ provided $y\ne 0$,
$z\not\in \{0,\pm \Ii\}$.  Irreducibility requires $z\not\in\{\pm \zeta,\pm \overline{\zeta}\}$.
Irreducible $\rho$ with $zz'=-1$ are  isomorphic. Here, $\alpha_1=\beta_1=\beta_2=1$,
so Tr$\,\FB=-\frac{3}{2}$ and $\FB=\textrm{diag}(t,-\frac{3}{2}-t)$
for $z=e^{2\pi \Ii t}$ and some $t\not\equiv\pm\frac{1}{4}$ (mod 1). Then as before
\begin{eqnarray}
\mathcal{X}&=&\left(\begin{array}{cc} 24\frac{20t^2+51t+32}{4t+3}&\frac{384x(3t+2)
(3t+1)(6t+5)(6t+7)}{(4t+5)(4t+3)^2}\\
\frac{384}{x(4t+1)}&-12\frac{40t^2+18t+1}{4t+3}\end{array}\right)\,,\\
\label{eq:fund2b}
{\Xi}(\FC)&=&(1-\FC)^{1/3}\left(\begin{array}{cc} 
(-1728\FC)^{-t}f(t+\frac{1}{3};\frac{5}{6};\FC)&
\mathcal{X}_{12}f(t+\frac{4}{3};\frac{5}{6};\FC)\\ 
\mathcal{X}_{21}f(-\frac{1}{6}-t;\frac{5}{6};\FC)&
f(-t-\frac{7}{6};\frac{5}{6};\FC)\end{array}\right)\,,\\
y&=&x\frac{\sqrt{3}}{6912\,\,\,432^{2t}}\frac{\Gamma(2t+\frac{3}{2})^2}{\Gamma(2t+\frac{4}{3})\,\Gamma(2t+\frac{2}{3})}\,.
\end{eqnarray}
Again, for irreducible $\rho$, any possible $t$ yields bijective $\FB$.
The 4 indecomposable but reducible $\rho$ are:

\smallskip\noindent$z=\zeta$: here, $\rho=\varmu_{2}\sdsum\varmu_{4}$ and $\FB$
is bijective iff $t\in\frac{1}{6}+\bbZ_{\ge 0}$;

\smallskip\noindent$z=-{\zeta}$: here $\rho=\varmu_{8}\sdsum\varmu_{10}$ 
and $\FB$ is bijective iff $t\in\frac{2}{3}+\bbZ_{\ge 0}$;

\smallskip\noindent$z=\overline{\zeta}$: here, $\rho=\varmu_{10}\sdsum\varmu_{8}$
and $\FB$ is bijective iff $t\in-\frac{1}{6}+\bbZ_{\ge 0}$;

\smallskip\noindent$z=-\overline{\zeta}$: here $\rho=\varmu_{4}\sdsum\varmu_{2}$ and $\FB$
is bijective iff $t\in\frac{1}{3}+\bbZ_{\ge 0}$.

\noindent Again the `missing' $t$ correspond to other reducible $\rho$.

The holomorphic analysis is identical to before. Here Tr$\,\lambda^{hol}=\frac{1}{2}$ or
$\frac{3}{2}$, depending on whether or not Re$\,t\le\frac{1}{2}$. In the former case $w^{(i)}=(2,4)$,
and in the latter it equals (8,10). $\bbX^{(i)}(\tau)$ is as before. The indecomposable $\rho$ behave
exactly as before.

\smallskip
The final class of two-dimensional representations has det$\,T=\overline{\zeta}$:
\begin{equation}
T=\left(\begin{array}{cc}z&0\\ 0&\overline{\zeta} z^{-1}\end{array}\right)\ ,\qquad S=\left(
\begin{array}{cc}
\frac{\zeta z}{z^2-\overline{\zeta}}&y\\ \frac{(z^2-1)(z+\overline{\zeta})(z-\overline{\zeta})}{y
(z^2-\overline{\zeta})^2}&\frac{-\zeta z}{z^2-\overline{\zeta}}\end{array}\right)\,,
\end{equation}
for arbitrary $y,z\in\bbC$ provided $y\ne 0$,
$z\not\in \{0,\pm e^{-\pi \Ii/6}\}$. Irreducibility requires $z\not\in\{\pm 1,\pm
\overline{\zeta}\}$; irreducible $\rho$ related by $zz'=\overline{\zeta}$ are equivalent.
 Here, $\alpha_1=\beta_0=\beta_2=1$,
so Tr$\,\FB=-\frac{7}{6}$ and  $\FB=\textrm{diag}(t,-\frac{7}{6}-t)$
for some $t$ satisfying $z=e^{2\pi \Ii t}$. Then
\begin{eqnarray}
\mathcal{X}&=\left(\begin{array}{cc} 24\frac{t(60t+71)}{12t+7}&\frac{10368x
t(2t+1)(3t+2)(6t+7)}{(12t+13)(12t+7)^2}\\
\frac{10368}{x(12t+1)}&-4\frac{(6t+7)(60t-1)}{12t+7}\end{array}\right)\,,\\
\label{eq:fund2c}
{\Xi}(\FC)&\,=\left(\begin{array}{cc} f(t;\frac{7}{6};\FC)&
\mathcal{X}_{12}f(t+1;\frac{7}{6};\FC)\\ 
\mathcal{X}_{21}f(-\frac{1}{6}-t;\frac{7}{6};\FC)&
f(-t-\frac{7}{6};\frac{7}{6};\FC)\end{array}\right)\,,\\
y=&\,x\frac{\sqrt{3}\,2^{1/3}}{10368\,\,\,432^{2t}}\frac{\Gamma(2t+\frac{7}{6})^2}{\Gamma(2t)\,\Gamma(2t+\frac{4}{3})}\,.\end{eqnarray}
Again, for irreducible $\rho$, any possible $t$ yields bijective $\FB$.
The 4 indecomposable but reducible $\rho$ are:

\smallskip\noindent$z=1$: then $\rho=1\sdsum\varmu_{10}$ and $\FB$
is bijective iff $t\in 1+\bbZ_{\ge 0}$;

\smallskip\noindent$z=-1$: then $\rho=\varmu_{6}\sdsum\varmu_{4}$ 
and $\FB$ is bijective iff $t\in\frac{1}{2}+\bbZ_{\ge 0}$;

\smallskip\noindent$z=\overline{\zeta}$: then $\rho=\varmu_{10}\sdsum 1$
and $\FB$ is bijective iff $t\in-\frac{1}{6}+\bbZ_{\ge 0}$; 

\smallskip\noindent$z=-\overline{\zeta}$: then $\rho=\varmu_{4}\sdsum\varmu_{6}$ and $\FB$
is bijective iff $t\in\frac{1}{3}+\bbZ_{\ge 0}$.

\noindent Again the `missing' $t$ correspond to the other reducible
$\rho$.

The holomorphic story for irreducible $\rho$ is as before. Here Tr$\,\lambda^{hol}=\frac{5}{6}$
or $\frac{11}{6}$ depending on whether or not Re$\,t\le \frac{5}{6}$. This means $w^{(i)}$
will equal (4,6) or (10,12), respectively. The only new phenomenon here is the indecomposable
at $z=\overline{\zeta}$: at the end of Section 3.7 we learned that dim$\,\cM^{\lambda^{hol}}_0(
\varmu_{10}\sdsum 1)$ is 0, not 1. We find $w^{(i)}=(4,6)$.

\smallskip
The $\lambda^{hol}$-holomorphic two-dimensional theory is also studied in \cite{bib:M2}, though without quantifying the
relation between Fourier coefficients and the matrix $S$ (i.e. his statements are only basis-independent),
which as we see involves the Gamma function. Our two-dimensional story can be
extended to any triangle group \cite{bib:jBG}.

\subsection{vvmf in dimensions $<6$} 

Trivially, the spaces of vvmf for an arbitrary $\rho$ are direct sums of those for its indecomposable
summands. Theorem 3.6 reduces understanding the vvmf for an indecomposable $\rho$,
to those of its irreducible constituents.
In this section we prove any admissible $(\rho,w)$ is \textit{tight}, provided $\rho$ is irreducible
and of dimension $<6$ (recall the definition of tight at the end of Section 3.4). This means in dimension $<6$ we get all kinds of things for free (see Proposition 3.3), including identifying
the Hilbert--Poincar\'e series $H^\lambda(x;\rho)$.  These series were first computed in \cite{bib:Mar}, for the special case
of exponent $\lambda=\lambda^{hol}$, when $T$ is unitary, and the representation $\rho$ is what Marks calls
\textit{T-determined}, which means that any indecomposable $\rho'$ with the same $T$-matrix
is isomorphic to $\rho$. It turns out that most  $\rho$ are $T$-determined. We will see this
hypothesis is unnecessary, and we can recover and generalise his results with much less effort.
The key observation is the following, which is of independent interest:

\medskip\noindent\textbf{Theorem 4.1.} \textit{Let $(\rho,w)$ be admissible and $T$ diagonal.
Assume $\rho$ is irreducible and the dimension $d<6$. Then $\rho$ is} tight: \textit{an
exponent $\lambda$ is bijective for $(\rho,w+2k)$ iff Tr$\,\lambda=c_{(\rho,w+2k)}$.}\medskip

\noindent\textit{Proof.} The case $d=1$ is trivial, and $d=2$ is explicit in Section 4.2, so it suffices
to consider $d=3,4,5$. Without loss of generality (by tensoring with $\varmu_{-w}$)
we may assume $\rho$ is a true representation of $\overline{\Gamma}$, i.e. that $(\rho,0)$
is admissible. Let $\alpha_i=\alpha_i(\rho,0)$, $\beta_j=\beta_j(\rho,0)$. Let $\bbX^{(i)}(\tau)$
be the free generators which exist by Theorem 3.4(a), and let $w^{(1)}\le\cdots
\le w^{(d)}$ be their weights. We will have shown that $\rho$ is tight, if we can show that
these $w^{(i)}$ agree with those in the tight Hilbert--Poincar\'e series \eqref{eq:fake}. This is 
because this would require all dim$\,\cM_{w+2k}^\lambda(\rho)$ to equal that predicted
by $H_{tt}(x)$, as given in Proposition 3.3, and that says $\lambda$ will be bijective
iff $\lambda$ has the correct trace. In fact it suffices to verify that the values of $w^{(i)}-w^{(1)}$
match the numerator of \eqref{eq:fake}, as the value of $\sum_iw^{(i)}$ would then also
fix $w^{(1)}$.

Let $n_i$ be the total number of generators $\bbX^{(i)}(\tau)$ with weight $w^{(i)}\equiv 2i$ (mod 12).
By Theorem 3.4(b) we know $\sum_in_i=d$, $n_i+n_{2+i}+n_{4+i}=\alpha_i$, $n_j+n_{3+j}=\beta_j$ for all $i,j$. These have solutions
\begin{equation}\label{lineqsoln}
(n_0,n_1,n_2,n_3,n_4,n_5)=(\alpha_0-\beta_2+t-s,\beta_1-s,\beta_2-t,s-t+\alpha_1-\beta_1,s,t)
\end{equation}
for parameters $s,t$.

Consider first $d=3$. The inequalities \eqref{eq:inequal} force $\beta_i=1$ and $\{\alpha_0,
\alpha_1\}=\{1,2\}$. By Proposition 3.3(a), we can assume without loss of generality (hitting
with $\varmu_{6}$ if necessary) that $\alpha_1=2$. Then the only nonnegative solutions
to \eqref{lineqsoln} are $(n_i)=(0,1,1,1,0,0),(1,1,0,0,0,1),(0,0,0,1,1,1)$. From \eqref{trmod6}
we see $L:=\mathrm{Tr}\,\lambda\in\bbZ$. Theorem 3.4(b) says $\sum_iw^{(i)}=12L$. The
inequality \eqref{eq:bounds} requires $w^{(1)}=4L-2$, so the only possibility consistent
with the given values of $n_i$ together with $\sum_iw^{(i)}=12L$ is $(w^{(i)})=
(4L-2,4L,4L+2)$, which is the prediction of \eqref{eq:fake}.

Consider next $d=4$. As before we may assume $(\beta_i)=(2,1,1)$ and $(\alpha_i)=(2,2)$.
Then \eqref{lineqsoln} forces $(n_i)=(1,1,1,1,0,0),(1,0,0,1,1,1,1),(2,1,0,0,0,1),(0,0,1,2,1,0)$.
Again $L:=\mathrm{Tr}\,\lambda\in\bbZ$ and $\sum_iw^{(i)}=12L$, and  \eqref{eq:bounds}
forces $w^{(1)}\in\{3L-3,3L-1\}$ (if $L$ is odd) or $w^{(1)}\in\{3L-2,3L\}$ (if $L$ is even).

We claim for each $L$ there is a unique possibility for the $w^{(i)}$ which is compatible
with $\sum_iw^{(i)}=12L$, the listed possibilities for $(n_i)$, the 2 possible values for
$w^{(1)}$ given above, and the absence of a `gap' in the sense of \eqref{eq:gap}.
When $L$ is even this is $(w^{(i)})=(3L-2,3L,3L,3L+2)$; when $L$ is odd this is $(w^{(i)})=(3L-3,3L-1,3L+1,3L+3)$. These match \eqref{eq:fake}.

Finally, consider  $d=5$. We may take $(\beta_i)=(1,2,2)$ and $(\alpha_i)=(3,2)$,
so 
\begin{equation}(n_i)\in\{(1,2,2,0,0,0),(1,1,1,0,1,1),(1,0,0,0,2,2),(0,1,2,1,1,0),(0,0,1,1,2,1)\}\,,\end{equation}
$L:=\mathrm{Tr}\,\lambda\in\bbZ$ and $\sum_iw^{(i)}=12L$.  

If $L=5L'$, then \eqref{eq:bounds} forces $w^{(1)}=12L'-4$ or $12L'-2$. We find the
only possible value of $w^{(i)}$ is $(12L'-4,12L'-2,12L',12L'+2,12L'+4)$.

If $L=5L'+1$, then \eqref{eq:bounds} forces $w^{(1)}=12L'$. We find the
only possible value of $w^{(i)}$ is $(12L',12L'+2,12L'+2,12L'+4,12L'+4)$.

If $L=5L'+2$, then \eqref{eq:bounds} forces $w^{(1)}=12L'+2$ or $12L'+4$. We find the
only possible value of $w^{(i)}$ is $(12L'+2,12L'+4,12L'+4,12L'+6,12L'+8)$.

If $L=5L'+3$, then \eqref{eq:bounds} forces $w^{(1)}=12L'+4$ or $12L'+6$. We find the
only possible value of $w^{(i)}$ is $(12L'+4,12L'+6,12L'+8,12L'+8,12L'+10)$.

If $L=5L'+4$, then \eqref{eq:bounds} forces $w^{(1)}=12L'+6$ or $12L'+8$. We find the
only possible value of $w^{(i)}$ is $(12L'+8,12L'+8,12L'+10,12L'+10,12L'+12)$.

All of these agree with \eqref{eq:fake}.\qed\medskip

Proposition 3.3 gives some consequences of tightness.

\subsection{Further remarks}

Now let's turn to more general statements. The simplest way to change the weight $w$ has already
been alluded to several times. Namely, suppose we are given any admissible multiplier system $(\rho,w)$ with bijective
 $\FB$ and fundamental
matrix ${\Xi}(\tau)$. Recall  the multiplier $\varmu_w$ of
$\Delta^{w}(\tau)$. Then for any $w'\in\bbC$, $(\varmu_{w'}\otimes\rho,w+12w')$ is
admissible, with bijective  and fundamental matrices 
$\FB+{w'}1_d$  and $\Delta^{w'}(\tau)\,
{\Xi}(\tau)$.

Suppose bijective $\Lambda,\Lambda'$ with corresponding fundamental matrices $\Xi(\tau),\Xi'(\tau)$ are known for admissible $(\rho,w)$ and $(\rho',w')$. Then the $dd'$ columns
of the Kronecker matrix product $\Xi(\tau)\otimes\Xi'(\tau)$ will manifestly lie in $\cM^!_{w+w'}(
\rho\otimes\rho')$, and will generate over $\bbC[J]$ a full rank submodule of it. By Proposition
3.2 the differential operators $\nabla_i$ then generate from that submodule all of
$\cM^!_{w+w'}(\rho\otimes\rho')$. In that way, bijective exponents and fundamental matrices
for tensor products (and their submodules) can be obtained.

The easiest and most important products involve  the six 
one-dimensional $\overline{\Gamma}$-representations $\varmu_{2i}$. Here we can be much more explicit.
Equivalently, we can describe the effect of changing the weights by even integers but
keeping the same representation. For simplicity we restrict to even integer weights.

\medskip\noindent\textbf{Proposition 4.1.} \textit{Let $(\rho,0)$ be admissible and $T$ diagonal.
Fix a bijective $\FB$, with corresponding ${\Xi}(\tau)$. Then for 
$i=2,3,4,5$ respectively, the columns for a fundamental matrix for $(\rho,2i)$ can be obtained 
as a linear combination over $\bbC$ of the columns of, respectively,:}

%\smallskip\noindent\textbf{(1)} \textit{$\Der{\Xi}=q^\Lambda(\Lambda+\cdots)$, 
%$E_4^2\left(E_4\Der{\Xi}-E_6{\Xi}\FB\right)/\Delta=q^\Lambda(1728\cA_2+\cdots)$,
%$(E_6^2\Der{\Xi}-E^2_4E_6{\Xi}\FB)/\Delta=q^\Lambda(-1728\cA_3+\cdots)$;}

\smallskip\noindent\textbf{(2)} \textit{$E_4{\Xi}=q^\Lambda(1_d+\cdots)$ and
$\Der^2{\Xi}-E_4{\Xi}\FB(\FB-\frac{1}{6})=q^\Lambda(1728\cA_3(\cA_3-\frac{1}{3})q+\cdots)$}\,;
  
\smallskip\noindent\textbf{(3)} \textit{$E_6{\Xi}=q^\Lambda(1_d+\cdots)$ and $E_4\Der{\Xi}-E_6{\Xi}\FB=
q^\Lambda(1728\cA_2 q+\cdots)$}\,;

\smallskip\noindent\textbf{(4)} \textit{$E_4^2{\Xi}=q^\Lambda(1_d+\cdots)$ and
$E_6\Der{\Xi}-E_4^2{\Xi}\FB=q^\Lambda(1728\cA_3 q+\cdots)$}\,;

\smallskip\noindent\textbf{(5)} \textit{$E_4E_6{\Xi}$, 
$E_4\left(E_4\Der{\Xi}-E_6{\Xi}\FB\right)=q^\Lambda(1728\cA_2q+\cdots)$, and
$E_6\left(\Der^2{\Xi}-E_4{\Xi}\FB(\FB-\frac{1}{6})
\right)=q^\Lambda(1728\cA_3(\cA_3-\frac{1}{3})q+\cdots)$}\,.

\medskip\noindent\textit{Proof.} Let $\alpha_i=\alpha_i(\rho,w)$ and $\beta_j=\beta_j(\rho,w)$,
and write $\FB_{(i)}$ for some to-be-determined bijective exponent at weight $w+2i$. Then
Tr$\,\FB_{(i)}-\mathrm{Tr}\,\FB$ can be read off from \eqref{tracek}.
%=\alpha_1-\beta_0,\beta_2,\alpha_1,\beta_1+
%\beta_2,\alpha_1+\beta_2$ for  $i=1,2,3,4,5$ respectively, where $\alpha_i,\beta_j$ are the 
%multiplicities for $\cA_{(w)},\cB_{(w)}$. 
Consider first the case $i=5$.  Theorem 3.3(c) says 
$\cA_3(\cA_3-\frac{1}{3})$ has rank $\beta_2$, while $\cA_2$ has rank
$\alpha_1$. The column spaces of the second and third matrices given above
for $i=5$ 
have trivial intersection, since any $\FA(\tau)$ in their intersection is a vvmf for $(\rho,w+2)$
which would be divisible by $E_4(\tau)^2E_6(\tau)$. This would mean 
$E_{14}(\tau)^{-1}\Delta(\tau)\FA(\tau)\in\cM^!_w(\rho)$ would lie in the kernel of $\FE_\FB$, so by bijectivity
$\bbX(\tau)$ must be 0. The desired generators will be linear combinations of $\alpha_1$
columns of the second matrix with $\beta_1$ of the third and $d-\beta_1-\alpha_1$ of the
first. These would define a matrix $\Xi_{(5)}(\tau)$ whose $\Lambda_{(5)}$ has the correct
trace, and therefore it must be a fundamental matrix (since its principal part map $\cP_{\Lambda_{(5)}}$ is
manifestly surjective).

The other cases listed are easier. \qed\medskip

In all these cases $2\le i\le 6$, Proposition 4.1 finds a bijective exponent $\Lambda_{(i)}$
such that $\Lambda_{(i)}-\Lambda$ consists of 0's and 1's (for $i=6$, $\Lambda_{(6)}=\Lambda
+1_d$ always works).

The case $i=1$ is slightly more subtle. Note that a fundamental matrix for $w+2i+12j$ is $\Delta
(\tau)^j$ times that for $w+2i$. So one way to do $i=1$ is to first find $i=4$, then find the
$i=3$ from $w+8$, then divide by $\Delta(\tau)$. Here is a more direct approach: almost always,
the columns of the fundamental matrix for $i=1$ is a linear combination over $\bbC$ of the
columns of $M_1(\tau):=\Der{\Xi}$, 
$M_2(\tau):=E_4^2\left(E_4\Der{\Xi}-E_6{\Xi}\FB\right)/\Delta=q^\Lambda(1728\cA_2+\cdots)$, and
$M_3(\tau):=(E_6^2\Der{\Xi}-E^2_4E_6{\Xi}\FB)/\Delta=q^\Lambda(-1728\cA_3+\cdots)$. Indeed, take
any $v\in\mathrm{Null}(\cA_2-\frac{1}{2})\cap
\mathrm{Null}\left((\cA_3-\frac{1}{3})(\cA_3-\frac{2}{3})\right)$. That intersection has dimension at least $\alpha_1-\beta_0$, because Null$(\cA_2-\frac{1}{2})$
has dimension $\alpha_1$ and  Null$(\cA_3-\frac{i}{3})$ has dimension $\beta_i$.
Then $M_2v=q^\Lambda(v/2+\cdots)$ and
$M_3v=q^\Lambda(jv/3+\cdots)$ where $j=1$ or 2, so $2jM_2(\tau)v-3M_3(\tau)v\in q^(\Lambda+1_d)\bbC^d[[q]]$. Generically, this gives $\alpha_1-\beta_0$ linearly independent vvmf, 
which together with $d-\alpha_1+\beta_0$ columns of $M_1(\tau)$ will give the columns of
the $i=1$ fundamental matrix. This method doesn't work for all $\rho$ however, e.g. it fails for $\rho=1$.

\begin{acknowledgement} I would like to thank Peter Bantay for a fruitful collaboration, Chris Marks,
 Geoff Mason, and Arturo Pianzola for conversations, and my students Jitendra Bajpai and
 Tim Graves.  My research is supported in part by NSERC. 
\end{acknowledgement}

%\bibliographystyle{plain}
%\bibliography{/home/bantay/bibfiles/CFT,/home/bantay/bibfiles/Moonshine,/home/bantay/bibfiles/Riemannsurfaces,/home/bantay/bibfiles/VOA,/home/bantay/bibfiles/conformalcharacters,\string"/home/bantay/bibfiles/modular functions\string",/home/bantay/bibfiles/number_theory,/home/bantay/bibfiles/modularrepresentation,/home/bantay/bibfiles/my}

\begin{thebibliography}{10}

\begin{footnotesize}

\bibitem{bib:AM} G. Anderson and G. Moore, ``Rationality in conformal field theory'',  Commun. 
Math. Phys.  117  (1988) 441--450. 


\bibitem{bib:AS}  {A.\ O.\ L.\ Atkin, H.\ P.\ F.\ Swinnerton-Dyer} 
{Modular forms on noncongruence subgroups.} In: Proc.\ Symp.\
Pure Math.\ {\bf 19} (AMS, Providence, 1971), ed.\ by T.\ S.\ Motzkin, pp.1--26.


\bibitem{bib:jBG} J. Bajpai and T. Gannon, in preparation.


\bibitem{bib:BG} P. B\'antay and T. Gannon, {}``Conformal characters
and the modular representation'', JHEP 0602 (2006) 005.

%\bibitem{bib:GMM} T. Gannon, C. Marks, and G. Mason, ``Cohomology of
%the modular group'', in preparation.

\bibitem{bib:BG2} P. B\'antay and T. Gannon, {}``Vector-valued modular
functions for the modular group and the hypergeometric equation'',
Commun. Number Th. Phys. 1 (2008) 637--666.

\bibitem{bib:BG3} P. B\'antay and T. Gannon, in preparation.


\bibitem{bib:Bolicm}  A. A. Bolibruch, ``The Riemann-Hilbert problem and Fuchsian
differential equations on the Riemann sphere'', Proc. Intern. Congr. Math. (Z\"urich, 1994),
Birkh\"auser, Basel 1995, pp.1159--1168. 

\bibitem{bib:Bolrev}  A. A. Bolibruch, ``The Riemann-Hilbert problem'', Russian Math. Surveys
45:2 (1990) 1--47.

\bibitem{bib:Borlift} R. E. Borcherds, ``Automorphic forms with singularities on Grassmannians'', 
 Invent. Math.  132  (1998) 491--562.

\bibitem{bib:Bor} R. E. Borcherds, {}``Gross-Kohnen-Zagier theorem
in higher dimensions'', Duke Math. J. 97 (1999) 219--233; correction
Duke Math. J. 105 (2000) 183--184.

\bibitem{bib:Za} J. H. Bruinier, G. van der Geer, G. Harder, and D. Zagier,
\textit{The 1-2-3 of Modular Forms} (Springer, Berlin, 2008).

\bibitem{bib:DMMV} R. Dijkgraaf, J. Maldecena, G. Moore, and E. Verlinde, ``A black
hole Farey tale'', preprint arXiv: hep-th/0005003.

\bibitem{bib:EZ} M. Eichler and D. Zagier, {}\textit{The Theory of Jacobi Forms},
 Prog. Math. 55 (Birkh\"auser, Boston, 1985).


\bibitem{bib:Ga} T. Gannon, in preparation.

\bibitem{bib:Hi} E. Hille, \textit{Ordinary Differential Equations in the
Complex Domain}
(Wiley \& Sons, 1976).


\bibitem{bib:Kim} T. Kimura, ``On the Riemann problem on Riemann surfaces'', Lecture Notes
in Math. 243 (Springer, 1971), pp.218--228.


\bibitem{bib:KM} M. Knopp and G. Mason, {}``Vector-valued modular
forms and Poincar\'e series'', Illinois J. Math. 48 (2004) 1345--1366.

\bibitem{bib:KMlog} M. Knopp and G. Mason, {}``Logarithmic vector-valued modular forms'', 
Acta Arith. 147 (2011) 261--282; arXiv:0910.3976.

\bibitem{bib:KrMa} M. Krauel and G. Mason, 
``Vertex operator algebras and weak Jacobi forms'',
Internat. J. Math. 23 (2012), no. 6, 1250024. 

\bibitem{bib:Lev} A. H. M. Levelt, ``Hypergeometric functions II', Indag. Math. 23 (1961) 373--385.

\bibitem{bib:Mar} C. Marks,  ``Irreducible vector-valued modular forms of dimension less than six,'' 
Illinois J. Math.  55 (2011) 1267--1297; arXiv:1004.3019.

\bibitem{bib:MM} C. Marks and G. Mason, ``Structure of the module of vector-valued
modular forms'', J. London Math. Soc. 82 (2010) 32--48; arXiv:0901.4367. 

\bibitem{bib:M2} G. Mason, ``2-dimensional vector-valued modular forms'', Ramanujan J. 17
(2008) 405--427.

\bibitem{bib:Mas} G. Mason, ``Vector-valued modular forms and linear
differential operators'', Intl J. Number Th. 3 (2007) 377--390.

\bibitem{bib:MMS} S. D. Mathur, S. Mukhi and A. Sen,  ``Differential equations for correlators and characters in arbitrary rational conformal field theories,''
Nucl. Phys. B312 (1989) 15--57. 

\bibitem{bib:Mil} A. Milas, ``Virasoro algebra, Dedekind eta-function and specialized Macdonald's 
identities,''  Transform. Groups  9 (2004) 273--288.

\bibitem{bib:miy} T. Miyaki, \textit{Modular Forms} (Springer, Berlin Heidelberg,
2006).

\bibitem{bib:Miy1} M. Miyamoto, ``A modular invariance on the theta functions defined on vertex operator algebras'', Duke Math. J. 101 (2000) 221--236.

\bibitem{bib:Miy} M. Miyamoto, {}``Modular invariance of vertex
operator algebras satisfying $C_{2}$-cofiniteness'', Duke Math.
J. 122 (2004) 51--91.

%\bibitem{bib:Pet} H. Petersson, ``Neuere Untersuchungen \"uber automorphe Formen
%komplexer Dimension. Bericht'', Jahr. der Deutschen Math.-Verein.


\bibitem{bib:Roh} H. R\"ohrl, ``Das Riemann-Hilbertsche Problem der Theorie der linearen
Differentialgleichungen'', Math. Ann. 133 (1957) 1--25.

\bibitem{bib:Se} A. Selberg, ``On the estimation of Fourier coefficients
of modular forms'', In: Theory of Numbers. \textit{Proc. Sympos. Pure Math.} VIII (Amer. Math. Soc.,
Providence, 1965) pp.1--15.

\bibitem{bib:Shim} G. Shimura, \textit{Introduction to the Arithmetic
Theory of Automorphic Functions} (Princeton University Press, 1971).


\bibitem{bib:Sk}  N.-P. Skoruppa, {}``\"Uber den Zusammenhang zwischen Jacobiformen
und Modulformen halbganzen Gewichts'', PhD Thesis, Universit\"at Bonn
(1984).

\bibitem{bib:TW} I. Tuba and H. Wenzl, {}``Representations of the
braid group $B_{3}$ and of SL($2,\mathbb{Z}$)'', Pacific J. Math.
197 (2001) 491--510.



\bibitem{bib:West} B. W. Westbury, ``On the character varieties of free products of cyclic groups'', 2001 preprint.

\bibitem{bib:Zh} Y. Zhu, {}``Modular invariance of characters of
vertex operator algebras'', J. Amer. Math. Soc. 9 (1996) 237--302.

\end{footnotesize}\end{thebibliography}

%\appendix
%dummy comment inserted by tex2lyx to ensure that this paragraph is not empty
%dummy comment inserted by tex2lyx to ensure that this paragraph is not empty
%dummy comment inserted by tex2lyx to ensure that this paragraph is not empty
%dummy comment inserted by tex2lyx to ensure that this paragraph is not empty
%dummy comment inserted by tex2lyx to ensure that this paragraph is not empty

\end{document}